\documentclass[11pt]{article}

\usepackage{import,mathrsfs,float,pgfplots,mhchem,
bbm,graphicx,subcaption,listings,xcolor,tcolorbox}\usepackage{etoolbox,booktabs}

\usepackage{caption}
\usepackage[acronym,toc,sort=use,nonumberlist]{glossaries}

\preto\printglossary{\glsaddall}
\makeglossaries

\lstdefinestyle{mathematica}{
    language=Mathematica,
    basicstyle=\small\ttfamily,
    backgroundcolor=\color{gray!10},
    frame=single,
    breaklines=true,
    commentstyle=\color{green!60!black},
    keywordstyle=\color{blue},
    stringstyle=\color{red},
    showstringspaces=false,
    tabsize=2
}
\usepackage[final,backref]{hyperref}
\usepackage{tikz-cd}
\usepackage{lineno}
\tikzcdset{every label/.append style = {font = \small}}
\usepackage{tikz}
\usetikzlibrary{graphs,decorations.pathmorphing,decorations.markings,arrows.meta,
automata,positioning}
\usepackage{subcaption}
\usepackage{nomencl}
\makenomenclature

\usepackage{makecell}
\usepackage{amssymb, amsfonts,mathtools}
\usepackage{eurosym,graphicx,latexsym, amsmath,amsthm,times,nccmath,url}
\usepackage{listings,verbatim,cases,tabularx}
\usepackage{graphicx,float,fancybox}
\usepackage{epsfig,epstopdf,tikz}
\usepackage[sans]{dsfont}
\usetikzlibrary{positioning,arrows,arrows.meta,calc}

\newtheorem{theorem}{Theorem}
\newtheorem{definition}{Definition}
\newtheorem{open}{Problem}
\newtheorem{example}{Example}

\newcounter{tabcounter}
\newcounter{defcounter}
\newcounter{opcounter}
\newcounter{thmcounter}
\newcounter{excounter}

\newglossarystyle{unifiedstyle}{
  \setglossarystyle{list}
  
}

\makeatletter
\newcommand{\beXa}[1][]{%
  \def\@temparg{#1}%
  \ifx\@temparg\@empty
    \begin{example}%
  \else
    \refstepcounter{excounter}
    \begin{example}[#1]%
    \label{e:#1}
    \newglossaryentry{e\theexcounter}{%
      name={#1},%
      description={p.~\pageref{e:#1}},
      user1={\theexcounter},%
      user2={Example}
    }%
  \fi
}
\newcommand{\eeXa}{\end{example}}

\newcommand{\beD}[1][]{%
  \def\@temparg{#1}%
  \ifx\@temparg\@empty
    \begin{definition}%
  \else
    \refstepcounter{defcounter}%
    \begin{definition}[#1]%
    \label{d:#1}%
    \newglossaryentry{d\thedefcounter}{%
      name={#1},%
      description={p.~\pageref{d:#1}},
      user1={\thedefcounter},%
      user2={Definition}%
    }%
  \fi
}
\newcommand{\eeD}{\end{definition}}

\newcommand{\beT}[1][]{%
  \def\@temparg{#1}%
  \ifx\@temparg\@empty
    \begin{theorem}%
  \else
    \refstepcounter{thmcounter}%
    \begin{theorem}[#1]%
    \label{t:#1}%
    \newglossaryentry{t\thethmcounter}{%
      name={#1},%
      description={p.~\pageref{t:#1}},%
      user1={\thethmcounter},%
      user2={Theorem}%
    }%
  \fi
}
\newcommand{\eeT}{\end{theorem}}
\newcommand{\beO}[1][]{%
  \def\@temparg{#1}%
  \ifx\@temparg\@empty
    \begin{open}%
  \else
    \refstepcounter{opcounter}%
    \begin{open}[#1]%
    \label{op:#1}%
    \newglossaryentry{o\theopcounter}{%
      name={#1},%
      description={p.~\pageref{op:#1}},
      user1={\theopcounter},%
      user2={Open Problem}%
    }%
  \fi
}
\newcommand{\eeO}{\end{open}}

\newcommand{\betA}[1][]{%
  \def\@temparg{#1}%
  \ifx\@temparg\@empty
    \begin{table}[htbp]%
  \else
    \refstepcounter{tabcounter}%
    \begin{table}[htbp]%

    \label{A:\thetabcounter}%

    \newglossaryentry{A\thetabcounter}{%
      name={#1},
      description={p.~\pageref{A:\thetabcounter}},%
      user1={\thetabcounter},%
      user2={Table}%
    }%
  \fi
}

\newcommand{\eetA}{\end{table}}
\makeatother
\newcommand{\printboth}{%
  \glsaddall  
  \printglossary[style=unifiedstyle, title={Definitions, Theorems, Examples and  Problems}]%
}

\newtheorem{lemma}{Lemma}
\newtheorem{proposition}{Proposition}
\newtheorem{corollary}{Corollary}

\newtheorem{remark}{Remark}
\newtheorem{question}{Question}
\newtheorem{assumption}{Assumption}

\def\beL{\begin{lemma}}\def\eeL{\end{lemma}}
\def\beP{\begin{proposition}}\def\eeP{\end{proposition}}
\def\beC{\begin{corollary}}\def\eeC{\end{corollary}}
\def\beR{\begin{remark}}\def\eeR{\end{remark}}
\def\beQ{\begin{question}}\def\eeQ{\end{question}}
\def\beA{\begin{assumption}}\def\eeA{\end{assumption}}

\definecolor{funccolor}{RGB}{25,25,112}
\definecolor{desccolor}{RGB}{64,64,64}

\def\bep{\begin{pmatrix}}\def\eep{\end{pmatrix}}
\def\bev{\begin{vmatrix}}\def\eev{\end{vmatrix}}
\def\bea{\begin{eqnarray*}}\def\eea{\end{eqnarray*}}
\def\bc{\begin{cases}}\def\ec{\end{cases}}
\def\BEN{\begin{enumerate}}\def\EEN{\end{enumerate}}
\def\BI{\begin{itemize}}\def\EI{\end{itemize}}

\newcommand{\be}[1]{\begin{equation}\label{#1}}
\newcommand{\ee}{\end{equation}}
\newcommand{\beq}{\begin{eqnarray}}
\def\eeq{\end{eqnarray}}
\def\eqr{\eqref}\def\fr{\frac}\def\lbl{\label}
\def\Lra{\Longrightarrow}\def\cNGM{\cite{Diek,Van,Van08}}
\def\Eq{\Leftrightarrow}\def\b{\beta}
\newcommand{\E}{\ensuremath{\mathbb{E}}}

\newcommand{\R}{\mathbb{R}}  

\def\al{\alpha}

\def\f{\varphi}  
\def\ga{\gamma}\def\Ga{\Gamma}\def\G{\Gamma}
\def\la{\lambda}\def\La{\Lambda}
\def\si{\sigma}
\def\th{\theta}

\def\T{\widetilde}

\newcommand{\bff}[1]{{\mbox{\boldmath$#1$}}}

\def\x{\boldsymbol{x}}
\newcommand\y{\boldsymbol{y}}

\def\xy{\x-\y}

\def\v1{\vec {\bff 1}}

\def\mR{{\mathcal R}}

\providecommand{\pr}[1]{\left(#1\right)}
\providecommand{\pp}[1]{\left[#1\right]}

\newcommand*{\Scale}[2][4]{\scalebox{#1}{$#2$}}

\newcommand{\figu}[3]{%
\begin{figure}[H]
\centering
\includegraphics[scale=#3]{#1}
\caption{#2\label{f:#1}}
\end{figure}
}

\long\def\symbolfootnote[#1]#2{%
\begingroup
\def\thefootnote{\fnsymbol{footnote}}\footnote[#1]{#2}%
\endgroup}

\def\Rto{R_{{1}}^{\tilde{2}}}
\def\Rot{R_{{2}}^{\tilde{1}}}

\def\Rtc{R_{2,c}}

\def\sd{s_{0}}

\def\sE{s_{ee}}

\def\adm{admissible}
\def\and{antisymmetric}

\def\ch{characteristic polynomial}

\def\coe{coefficient}

\def\con{condition}
\def\coo{coordinate}

\def\dec{decomposition}

\def\eig{eigenvalue}\def\eval{eigenvalue}

\def\equ{equation}

\def\fno{from now on}
\def\for{\forall}

\def\Fr{Furthermore, }

\def\How{However, }
\def\ie{i.e. }
\def\im{\item}

\def\mbw{may be written as}

\def\nne{non-negative}

\def\oth{otherwise}

\def\para{parameter}

\def\Prf{{\bf  Proof:}}

\def\qu{\quad}
\def\QED{\hfill {$\square$}\goodbreak \medskip}
\def\rep{representation}
\def\resp{respectively}
\def\satd{satisfied}\def\satg{satisfying}\def\sats{satisfies}\def\saty{satisfy}
\def\sec{\section}\def\ssec{\subsection}\def\sssec{\subsubsection}

\def\var{variable}
\def\wk{well-known}
\def\wlo{w.l.o.g.}

\newcommand\CRN{chemical reaction networks}

\newcommand\RN{reaction network}

\def\sm{stoichiometric matrix}
\def\aut{autocatalytic}

\def\brn{basic reproduction number}
\def\com{compartment}
\def\DFE{disease free equilibrium}

\def\ME{mathematical epidemiology}

\def\mSt{multi-strain model}

\def\NGM{next generation matrix}

\def\repF{reproduction functions}

\def\jin{Jacobian matrix of the invasion vector field with respect to the invasion variables}
\def\jrn{Jacobian matrix of the resident vector field with respect to the resident variables}

\def\CEP{competitive exclusion principle}

\def\fp{fixed point}\def\bfp{boundary fixed point}\def\fbp{fixed boundary point}

\def\MH{Metzler-Hurwitz}\def\RH{Routh-Hurwitz}

\def\Lf{Lyapunov function}

\def\gLV{generalized Lotka-Volterra model}

\def\LCP{\text{linear complementarity problem}}
\def\LV{Lotka-Volterra}\def\LVK{Lotka-Volterra-Kolmogorov}

\def\regS{regular splitting}

\def\num{numerator}
 \def\DSR{directed species reaction graph} 

\newtcolorbox{keyresult}[1][]{
    colback=blue!5!white,
    colframe=blue!75!black,
    title=#1,
    fonttitle=\bfseries
}

\newtcolorbox{examplebox}[1][]{
    colback=green!5!white,
    colframe=green!75!black,
    title=#1,
    fonttitle=\bfseries
}

\newtcolorbox{proofbox}[1][]{
    colback=yellow!5!white,
    colframe=orange!75!black,
    title=#1,
    fonttitle=\bfseries
}

\title{On the Similarity between Epidemiologic Strains, Minimal Self-Replicable Siphons,  and autocatalytic cores
in (Chemical) Reaction Networks: 
Towards a Unifying Framework
}

\author{ Florin Avram$^{1}$,  Rim Adenane$^{2}$, Lasko Basnarkov$^{3}$, Andras Horvath$^4$}
\begin{document}

\maketitle
\begin{center}	

$^{1}$ Laboratoire de Math\'{e}matiques Appliqu\'{e}es, Universit\'{e} de Pau, Pau, France; avramf3@gmail.com \\

$^{2}$ Laboratoire d'Analyse, G\'{e}om\'{e}trie et Applications, d\'epartement des Math\'ematiques, Universit\'e Ibn-Tofail,  Kenitra, Maroc; rim.adenane@uit.ac.ma \\

 $^{3}$ Faculty of Computer Science and Engineering, SS Cyril and Methodius University, and Macedonian Academy of Sciences and Arts, Skopje, Macedonia ; lasko.basnarkov@finki.ukim.mk \\

$^{4}$ Department of Computer Science, University of Turin, Italy; horvath@di.unito.it

\end{center}

\begin{abstract}
{\bf Motivation}: We aim to study boundary stability and persistence of
   positive odes in mathematical epidemiology models by importing structural tools from
  \CRN.\\
This is largely a review work, which attempts to bring closer together the fields of mathematical epidemiology (ME), and chemical reaction networks (CRNs), based on several observations. We started by observing that  \emph{epidemiologic strains}, defined as disjoint blocks in either the Jacobian on the infected variables, or as blocks in the next generating matrix (NGM),   coincide,  in most of the examples we studied,  with either the set of    \emph{critical minimal siphons} or with  the set of \emph{minimal autocatalytic sets (cores)} in an underlying CRN (an exception, not fully understood yet, is offered in Section 5). We leveraged this to provide a  definition of the \DFE\ (DFE) face/infected set as the union of either all minimal siphons, or of all cores (they coincide always in our examples). Next, we provide  a proposed definition of ME models, as models which have a unique boundary fixed point on the DFE face, and for which the Jacobian of the infected subnetwork  admits a regular splitting, which allows defining the famous next generating matrix. We then define the \emph{ interaction graph on minimal siphons} (IGMS), whose vertices are minimal siphons, and whose edges indicate the existence of reactions producing species in one siphon from species in another. When this graph is acyclic, we say the model exhibits a \emph{Acyclic Minimal Siphon Decomposition} (AMSD). For AMSD models whose minimal siphons partition the infection species,   we show that the NGM   is block triangular after permutation, which implies the  classical max structure of the reproduction number $R_0$ for multi-strain  models.
In conclusion, using irreversible reaction networks, minimal siphons and  acyclic siphon decompositions we provide a natural bridge from CRN to  ME.\\
 We implement  algorithms to compute IGMS and detect AMSD in the \texttt{Epid-CRN} Mathematica package (\url{https://github.com/florinav/EpidCRNmodels}) (which contains already modules to identify minimal siphons, criticality, drainability, self-replicability, etc). Finally, we illustrate on several multi-strain ME examples   how the block structure induced by AMSD, and  the ME reproduction functions (as introduced in \cite{AAN,AAHK}, these are the eigenvalues of the NGM which are always positive, for any values of the non-infection variables),  allow expressing  boundary stability and persistence conditions by comparing growth numbers to $1$, as customary in ME. Note that while not addressing the general Persistence Conjecture mentioned in the title, our work  provides a systematic method for  deriving boundary instability conditions for a significant class of structured  models.
\end{abstract}


\tableofcontents
\printboth
\section{Introduction}
Population dynamics, ecology,  mathematical epidemiology (ME), virology, the theory of chemical reaction networks (CRN), to name just some biological interaction networks (BIN) subfields,  study all positive dynamical systems, and have similar preoccupations: the existence and multiplicity of equilibria, their local and global stability, the occurrence of bifurcations,   persistence, permanence, extinction, etc.

CRN offers a versatile mathematical framework for modelling positive systems in a
variety of disciplines; examples include
biochemical, pharmacokinetic, epidemiological, ecological, social, neural,
and multi-agent networks  \cite{Goutsias}.   Of special relevance to us is that CRN theory (CRNT)  investigates    mass–action systems,  which include ME models, and revealed  the roles of Petri graphs and siphons  in determining boundary behavior. For some recent papers combining both new CRNT results (including for the SIR  \ME\ (ME) model) and interesting historical references, see \cite{Plesa,GagPoly,XuGAS,XuNon,BloStaVas}.Let us mention also that CRNT is unavoidable  in the study of ODEs with many variables (tens, hundreds). Note however that a considerable part of  CRNT focuses on systems which enjoy various degrees of ``reversibility", which is absent in ME.

Ecology studies often Lotka–Volterra–Kolmogorov/generalized Lotka–Volterra (LVK/GLV)
 systems, whose behaviour might seem simpler than those in ME, due to the  linear complementarity principle (LCP); however, as shown by Smale, generalized Lotka–Volterra systems may have  already very complicated dynamical behaviors, and display chaos  \cite{Smale}. For some recent developments, including some  which tie Smale's work to CRNs and to ecology, see \cite{vano2006chaos,sprott2010elegant,wang2023lower,Susits,banaji2024bifurcations,
 chen2024stability,jafari2025simplest,
 Plesa}.

 One subfield in which CRNT and ecology come together is that of Chemical organisation theory \cite{dittrich2007chemical,heylighen2015chemical,veloz2017reaction,
 veloz2020complexity,heylighen2024chemical}.
These authors study both ODE and PDE models, and  touched also briefly the topic
 of ME models \cite{peter2020structure}.

Epidemiology is especially relevant to our unification quest due to the intricate
boundary behaviors encountered in the study of its multi-strain models, whose stability and persistence may sometimes be expressed via remarkably simple threshold conditions involving basic and  reproduction numbers and invasion numbers (these are essentially the particular case of the invasion rates used in ecology, obtained when the resident invariant set  of boundary faces reduce to single fixed points).

Ecology,  ME   and CRNT   share also common open problems, like   the notorious persistence and permanence conjectures, which they formulate however  differently.  Thus, CRNT considers only the case of weakly-reversible systems (see Section \ref{s:per} for a short review and further references). In ecology,  the accent is on general, non-polynomial rates, but possibly generating strongly monotone semi-flows \cite{ButlerMcGehee1979,Hutson,Schreiber,Thieme2003,HirschSmith,SmithThieme}. Finally,
in ME, which is the focus of this paper, the accent is on using  the often simpler reproduction and invasion numbers obtained via regular splitting (see Section \ref{s:reg}) and on  restricting to ME models, which need still to be defined  (note that both ME and  ecology models  lack a precise definition, in the literature preceding our papers).

The fact that the persistence and permanence conjectures
are only studied by each of the three fields  for restricted  structural classes suggests the interest of integrating their different perspectives; let us also
  note that CRNs have started been applied in ecology and related fields in a few recent papers \cite{veloz2021reaction,heylighen2022role,VAA,AAHK,
reyes2025functional}.

Our paper focuses brings together several key concepts from different fields:
\BEN \im Siphons and reaction networks from CRNT.
\im  The next generation matrix (NGM) method from ME, which is also related to the  regular splitting method in numerical analysis.
\im The use of the ME reproduction
functions -- see Definition \ref{d:reF}-- for expressing stability in terms or reproduction and invasion numbers/rates, where the second concept originated in ecology (see \cite{Schreiber, HofSch,feng2025theory,seno2025persistent} for  some references), but is different from the one used in ME, due to the application here of regular-splitting, for models where  resident invariant sets (see Definition \ref{d:prS}) of boundary faces reduce to single fixed points.
\EEN

  To argue our thesis that  not sufficient effort has  been invested  in cross-pollination (interdisciplinary integration) of CRNT, ecology and ME,  we add now some more information on the three points enumerated above.
\BEN \im The  reaction network (RN) representation of ODEs, and the concept of siphons
  are  essential in CRNT, and absent in ME, even though they seem intuitively important. In support of this, our main result below  started from the observation that there  is a correspondence between minimal siphons and strains in multi-strain epidemic models. This lead to the notion of
 interaction graph on  the minimal siphons (IGMS), and to Theorem \ref{t:par}: if IGMS is acyclic, and the minimal siphons do not intersect, then the NGM can be block-triangularized (recovering \wk\  results obtained previously only in examples).

 \im The  ``next generation matrix" (NGM) method (which expresses
 the stability of the ``disease free equilibrium'' (DFE) in terms of one inequality $R_0<1$, where $R_0$ is the spectral radius  of a certain matrix, obtained by ``regular splitting" of the Jacobian of the ``invading variables"), is omnipresent in ME, and  almost absent outside the field. But, we have shown
 in \cite{AAHJ,JA}, that it is equally efficient for CRN models.

 \im Our  final piece of evidence are multi-strain ME models. Note that while $R_0$ and NGM are typically presented as  closed chapters in ME, for the analog  stability problem of non-DFE boundary fixed points of multi-strain ME models,
  a general theory is definitely lacking. Our first observation here is   that in simple
  multi-strain  models, boundary fixed points are into a one to one correspondence with the minimal siphon facets ${S_i}, i=1,...,n$ on which they are situated, and we incorporate this observation in our definition of multi-strain models \ref{d:mtS}.  \Fr  it was observed in various multi-strain ME models that certain ``invasion numbers'' $R^{S_i}=R_{\T {S_i}^c}^{S_i}$, with definitions similar to  $R_0$,  play a central role in stability criteria, and we provide a definition for them in \ref{d:reF}.  Now minimal siphons are a \wk\ concept  in chemical reaction networks (CRN), but never used in ME and ecology (note also that invasion numbers are compared to $1$ in ME, and to $0$ in ecology).   We believe that for a general theory of multi-strain  models, in particular for progressing on problems like  whether k-strain models satisfy \LCP/\CEP\ (Open Problem \ref{o:LCP}), a better integration of CRN, ME  and ecology might be useful.  As a first step in this direction, we offer a general purpose Mathematica package Epid-CRN,
available at \url{https://github.com/florinav/EpidCRNmodels }, and illustrate its use throughout the paper.
\EEN

\paragraph{Contributions.}
\begin{itemize}
  \item We prove block-triangularity of the NGM when minimal siphons form an acyclic partition (Theorem~\ref{t:par}).
  \item We implement  detection of cycles of IGMS  in \texttt{Epid-CRN}, and illustrate the theory on multi-strain examples.
      \im Our paper may be viewed also as a users manual for EpidCRN,  and  also as an essay  around several not enough known  themes:

\BEN \im We outline first  an ``x–y  framework'' extending the NGM method,  and
connect it to the older theory of regular  splitting of matrices.

 \im  We suggest one possible rigorous
  definition of ME models,
see Definition \ref{d:ME}, by imposing the existence of a DFE, and   the existence of a regular splitting  of the \jin\  at the DFE.

\im  We  review the role played by reproduction functions \cite{AAN,AAHK} for the stability analysis of  multi-strain  ME models (these are basically
the reproduction  numbers,
in which the resident  variables $\y$ are left free).
\EEN
\end{itemize}

\paragraph{Organization.}
We start with some mathematical underpinnings in Section \ref{s:mat}, which include the regular splitting  of Metzler matrices, a fundamental result in ME.

 Section~\ref{s:IGMS}   introduces a new type of interaction graphs on minimal siphons (IGMS) and  provides a structural result for the NGM (Theorem~\ref{t:par}), which assumes the acyclicity of IGMS.

The first ME example in Section \ref{s:Rah}  revisits a ME model  of \cite{RahVac},
 which illustrates both the important complementarity principle of Lotka-Volterra-Kolmogorov models,   and revisits the reproduction functions
 introduced in \cite{AAN,AAHK}.

 Section \ref{s:IGMS} introduces the directed interaction graphs on minimal siphons (IGMS), and  establishes in Theorem \ref{t:par} that  next–generation matrix \(K=F V^{-1}\) inherits a block–triangular structure from $V$ whenever the minimal siphons form a partition of the DFE support-- see Theorem \ref{t:par}.

  Section~\ref{s:Gav}, following \cite{Agu07,GavishMusa,GavCEP,Gavnew}, revisits a quite  complex multi-strain example, with temporary cross-immunity,  with ADE and immunity waning, but without co-infection. Some stability results are easy, due to the    block-diagonal structure of the NGM,  but the domain of  Hopf bifurcations remains an open problem.

  Section \ref{s:Koz} presents a co-infection  ME model where the minimal siphons intersect; nevertheless, the NGM has a triangular structure. Furthermore, there are more reproduction functions than minimal siphons.   This  poses 
   challenges from the point of view of interpreting the stability/invasion results via  reproduction functions, which we leave for further work.

  Sections \ref{s:per} offers a brief review of persistence theory. We argue  here that its  open problems might be easier under the restriction to ME models (in particular when the resident invariant sets of all proper siphon faces reduce to single fixed points).   Also, we include a sketch of the theory of invasion graphs \cite{HofSch}, which we believe might throw a light on the open problems of persistence and permanence of the model in the preceding section.

Section \ref{s:ap}  offers further comments on the NGM method.

Finally, sections  \ref{s:rel}, \ref{s:SDAS}, \ref{s:DSR}  review further topics from CRNT
which we conjecture might become relevant for ME in the future. These might look at a first reading like vague similarities between disparate fields.  But, as as it happened to two of the authors \cite{VAA}, we discovered that certain CRN results    throw a whole new light on \ME. For example, it turns out  that the presence of certain submatrices in the stoichiometric matrix (which is independent of the choice of admissible rates) render automatic the presence of various types of instabilities \cite{Vas,VasSta,BloStaVas}, under ``parameter-rich kinetics" (which do not include mass-action, but do include ``Michaelis Menten saturated" kinetics).

\section{Mathematical underpinnings}\label{s:mat}

{
\ssec{Positive dynamical systems}

 \beD[Positive / non-negative dynamical system]
A dynamical system is called \emph{positive} \cite{rantzer2015scalable} or \emph{non-negative} \cite{haddad2010}
if the non-negative orthant
\[
\mathbb{R}^n_{\ge 0} := \{ x \in \mathbb{R}^n : x_i \ge 0,\; i = 1,\dots,n \}
\]
is forward invariant under the flow.
\eeD

  The BIN disciplines are all concerned with positive  dynamical systems, but neglect sometimes important results in their  ``sister sciences".

  In this essay, we are making the case that more sharing of information between the positive dynamical systems fields
   could be beneficial.
  \beR \label{r:pos} As a first argument, we ask the reader to compare the ODE description of our first example \eqr{ODERahv} to its ``reaction networks" (RN) representation \eqr{RNv}. The superiority
  of the latter, encapsulated in lines with biological meaning like $i_1+s \to 2~i_1$,
  seems clear to us:
   \BEN \im Firstly, reactions represent a clear phenomenological description: in the case above, ``one susceptible meets one infected, and two infected depart", is a precise description of the infection process.
   \im Note that after inputting the corresponding rates
  (rts), numerous packages will produce the ODE description (simply by ``adding up" the reactions), as well as various graphs illustrating the interactions in this network, numerical simulations, and, possibly in the future, also symbolic-numeric information on the fixed points and their stability (EpidCRN does that for simple examples).
  \im Using the RN representation, exploring the model by changing  the rates, or by adding/removing reactions can be achieved with minimal effort.
  \EEN

  For a second argument, we recall the well-known  CRNT result that assuming all free parameters are \nne, a positive polynomial ODE system may be recognized by the absence of  ``syntactically negative cross-terms" \cite{hun}. This result seems practically unknown outside CRNT, and has been reproved numerous times.
\eeR

  It seemed to us
   that one step towards advancing  unification of the BIN sciences could be via creating a common general purpose package (symbolic, numeric, simulation and estimation of parameters). Below,  we present the Mathematica package Epid-CRN,
available at \url{https://github.com/florinav/EpidCRNmodels }, which is used in solving all the examples of this paper.

\ssec{Notations}

We denote by $\mathbb R^n$, $\mathbb R_+^n$ and $\mathbb R_{\ge0}^n$ the
set of order $n$ column vectors of reals, positive reals and non-negative
reals, respectively. The usual euclidean norm is denoted by $\|
z\|_2^2$. The family $\{e_1,\ldots, e_n\}$ denotes the canonical basis of
the vector space $\mathbb R^n$. If $z \in \mathbb R^n$, we denote by $z_i$
the $i$th component of $z$. The notation $x\gg0$ means that $x \in\mathbb
R_+^n$.

The set of $n\times m$ matrices with real entries is denoted by $\mathbb
R^{n\times m}$. For a matrix $A\in\mathbb R^{n\times m}$ we denote by
$A(i,j)$ the entry at the row $i$, column $j$. $A^t$ denotes the transpose
of $A$, and $A^{-t}$ will denote the transpose of the inverse of $A$.  For
matrices $A,B$ we write $A \leq B$ if $A(i,j) \leq B(i,j)$ for all $i$ and
$j$, $A<B$ if $A \leq B$ and $A \neq B$, and $A \ll B$ if $A(i,j) < B(i,j)$
for all $i$ and $j$. Let $U\subset \mathbb R^{n\times m}$, then $\bar U$
denotes its closure with respect to the standard Euclidean topology. Let
$A\in\mathbb R^{n\times m}$ be a matrix, then $\ker(A)$ is its kernel or
null space, while $\Im A$ is its image space. Let $V:\mathcal X \to \mathbb
R$ be a function, then $\ker V:=\{x\in\mathcal X| V(x)=0\}$. Let $\mathcal
Y \subset \mathcal X$, then $\ker V|_{\mathcal Y}:=\{x\in\mathcal Y|
V(x)=0\}$. Let $x=(x_1,..,x_n)^t \in \mathbb R^n$, then
$\|x\|_\infty:=\max_i |x_i|$.  Let $A\in\mathbb R^{n\times n}$, then
$\sigma_i(A):= a_{ii}+\sum_{j\ne i} |a_{ij}|$, $i=1,..,n$, and
$\mu_\infty(A):=\max_i \sigma_i(A)$.

Let $S$ be a set, then $|S|$ denotes its cardinality.

\beD [spectral radius, spectral abscissa, and Hurwitz matrix]\cite{horn2012}
 A) The spectral radius of a matrix A  is defined by
$$\rho(A) = \max\{|\la| ,  \la \in  Sp(A)\},$$
where  $Sp(A)$ denotes the
spectrum of A.

B) The spectral abscissa, \ie\ the largest real part among the eigenvalues of $A$, is denoted by
$$s(A) = \max \left\{{\rm Re}(\lambda)\ : \ \lambda\in Sp(A)\right\}.$$

C) A matrix $A$ is said to be \textit{Hurwitz} if $s(A) < 0$.
\eeD

\beD [Metzler matrix]\cite{Berman,farina2011}  A Metzler matrix A is a matrix such that $i \neq j \Lra A(i, j) \ge 0.$  These
matrices are also called quasi-positive matrices. 
\eeD

\beL\label{lem:metzler} {\bf The ``half eigenvector" inequalities
}.
	Let $A$ be an irreducible Metzler matrix. Then, $s(A)$ is a simple eigenvalue of $A$ and there exists a unique (up to scalar multiple) vector $x > \vec  0_n$ such that $Ax = s(A) x$. Let $z \geq \vec  0_n$ be a given non-zero vector. \BEN \im
If $Az \leq \lambda z$ for some scalar $\lambda$, then $s(A) \leq \lambda$, with equality if and only if $Az = \lambda z$.
\im  If $Az \geq \lambda z$ and $Az \neq \lambda z$, for some scalar $\lambda$, then $s(A) > \lambda$.
\EEN
\eeL
The first half of the lemma is a direct consequence of the Perron--Frobenius Theorem for nonnegative matrices \cite
{Berman}. The second part can be obtained from a straightforward application of \cite[Theorem 2.1.11]{Berman}.

Scalar functions like the logarithm and quotients  will be applied to vectors componentwise:
$$\ln y=(\ln y_1,\ln y_2, ...), \fr {1}{y}=\left(\fr 1{y_1},\fr 1{y_2},...\right),
\fr {y^*}{y}=\left(\fr {y_{1}^*}{y_1},\fr {y_{2}^*}{y_2}, ...\right)$$


\subsection{Invasion criteria for  extinct coordinates via  sub--Jacobians
}

Let $\dot x=f(x)$ on $\mathbb{R}^n_{\ge 0}$ and let
$F_W=\{x_i=0:\ i\in W\}$ be a boundary face (the variables with indices $W$ are ``extinct" on $F_W$).
If $E\in F_W$ is a boundary equilibrium, write the Jacobian in block form
$Df(E)=\begin{bmatrix}J_{WW}(E)&*\\ *&*\end{bmatrix}$.
Define the \emph{invasion spectral abscissa}
\[
\alpha_W(E):=\max\{\,\Re\lambda:\ \lambda\in\sigma(J_{WW}(E))\,\}.
\]
Then:
\[
\Scale[.88]{\alpha_W(E)<0 \Rightarrow \text{attracting in the $W$–directions (no invasion)},\qquad
\alpha_W(E)>0 \Rightarrow \text{repelling in the $W$–directions (invasion)}.}
\]
For \textbf{LVK/GLV} systems (see Definition \ref{d:xyLV}) $\dot x_i=x_i\bigl(c_i+(Ax)_i\bigr)$ one has, on $F_W$,
\[
J_{WW}(E)=\mathrm{diag}\bigl(g_i(E)\bigr)_{i\in W},\quad
g_i(E):=c_i+(AE)_i,
\]
hence
\[
\alpha_W(E)=\max_{i\in W} g_i(E).
\]

{\bf Persistence criterion in a GLV.}
Because $\dot x_i$ has the factor $x_i$, each singleton $\{i\}$ is a minimal semilocking set. Persistence holds if every siphon face either has no boundary equilibrium or, when it does, satisfies $\alpha_W(E)>0$ (i.e. at least one extinct coordinate has positive per-capita growth at $E$).

\subsection{How to represent CRN and  ME  ODE models in a symbolic package?}
{For parameterizing ME and CRN ODE models},
   there are three natural choices:

 \BEN \im
 The parametrization   used traditionally for ODE models   is
 $X'=RHS(X)$. A model is thus  a pair $(X,RHS)$ where $X$ is the vector of the
 variables whose evolution we study and $RHS$ (right hand side) denotes the function that gives the derivatives.
 This parametrization is sufficient for studying small size symbolic bifurcation problems, but becomes awkward beyond that.

  \im In chemical reaction networks theory, $RHS$ is decomposed as \be{crnsys}X'=RHS(X)=\G~rts(X)\ee
where
  \BI  \im
   $\G$ is the ``stoichiometric matrix'' (SM), whose columns  represent directions in which several species/\com s change simultaneously,
   it is viewed as an embodiment of the ``structure" of the model.
   \im  $rts(X)$ is the  vector of rates of change associated to each direction (assumed all to be \nne), also known as kinetics. They are separated from the structure $\G$, since they are less certain; this fact, well accepted in CRNT,  is equally true in ME.
    \EI
    A model hence is defined as a triple $(X, \G, rts)$.

 \beXa[SIRS ODE without inflows and outflows (that is, without demography)] This is defined by the
 triple $(X, \G, rts)$:
  \be{SIRc}
X'=\bep
 s' \\  i'    \\  r'\eep
 =  \bep - 1&0&1&-1 \\
  1&- 1&0&0\\
  0&1&- 1&1
\eep \bep
\beta s  i    \\ \gamma_ i i    \\ \gamma_ r r \\ \gamma_ s s  \eep:=\Gamma ~rts(X).
\ee
\eeXa

Note that this representation is also a natural first step towards defining an associated CTMC model (continuous-time Markov chain) on the integers.

This parametrization has had remarkable successes, including recently for
studying \emph{robust questions} like the existence of \emph{robust Lyapunov
  functions} which depend only on the reaction rates
$rts(X)$~\cite{al2014robust,BG14,BlaGio17poly,Franco,AAScomp}. All these
papers assume admissibility conditions for the rates $rts(X)$ (see \cite{AdLS,ARS24}),
which make them behave essentially like polynomial mass-action systems.

\im {\bf The reaction network (RN) representation}. The third parametrization, via reactions, the traditional one used in CRNT, replaces
each column of $\G$ by a \emph{reaction/interaction}. For SIR, the reactions and rates are \resp:

\begin{center}
\begin{minipage}{0.45\linewidth}
\[
\mathrm{RN} = \left\{
\begin{aligned}
 S + I &\;\longrightarrow\; 2I,\\
 I &\;\longrightarrow\; R,\\
 R &\;\longrightarrow\; S,\\
 S &\;\longrightarrow\; R
\end{aligned}
\right\}.
\]
\end{minipage}%
\hfill
\begin{minipage}{0.45\linewidth}
\[
\mathrm{rts} = \left\{
\begin{aligned}
 &\beta S I,\\
 &\gamma_i I,\\
 &\gamma_r R,\\
 &\gamma_s S
\end{aligned}
\right\}.
\]
\end{minipage}
\end{center}

Note that  the first column in \eqr{SIRc} is replaced by $s + i \to 2 i$
(not by $s \to i$), which constitutes a precise phenomenological representation of
infections at the individual level: the meeting of a susceptible and an
infectious results in two infectious. This is called  mass-action representation.
\beD [mass-action reaction]\cite{Horn}.
A reaction is said to be mass-action if the coefficients multiplying the reactant species ($s + i$) in this case appear also as exponents in the rate of the reaction.
\eeD

The  other
columns in \eqr{SIRc},  represented like $i\to r$, etc, which correspond to transfers which do not occur after meeting individuals from other \com s, do not require explanation.

\beR \BI
\im  Essentially, the RN representation replaces $\G$ by the difference of two
\nne\  matrices (found on the LHS and RHS of each
reaction), such that \be{ab}\G=\beta-\alpha\ee
(which are called input/reactants and output/products matrices, \resp).

In fact, the RN representation may be viewed as an elegant way of visualizing
corresponding columns of $\alpha$ and $\beta,$ in parallel.

 \im Note that the RN representation \eqr{ab} is the first step towards
defining an associated discrete time Markov chain (DTMC) model on the
integers, where we distinguish between inflows and outflows in a state,
which is only possible in discrete time.

  It is surprising that even though the ODE only depends on $\G$, useful information about the ODE, like for example the siphons (see below) requires
  $\alpha$ and $\beta.$

 This gets even more surprising once we note that to one stoichiometric matrix $\G$ one may associate an
infinite family of $(\alpha,\beta)$ pairs which \saty\ $\beta-\alpha=\Ga$,
with different mass-action rates,  but the same ODE.

 \im The computation of siphons --see Definition \ref{d:sip}-- may be achieved either algorithmically, using the matrices $(\alpha,\beta)$, or, for small examples, by eye inspection.

\im The RN representation is  the one that best avoids
human errors,  when entering large systems with tens of
reactions.

\EI
\eeR
\EEN

In the RN \rep, each column of $\G$, that is, each \emph{reaction} is
associated to a directed pair consisting of a \emph{source complex}, and a
\emph{product complex}, which we proceed now to define.

\beD[CRN]\cite{feinberg2019} \label{d:CRN} A
CRN is defined by a triple $(\mathcal{S},\mathcal{C},\mathcal{R})$, where
$\mathcal{S},\mathcal{C},\mathcal{R}$ are the set of species, complexes and
reactions, respectively.

\begin{align*}
    \mathcal{S} &= \{s_1, ... , s_i,..., s_{|\mathcal{S}|}\}\\
    \mathcal{C} & = \{y_1,...,y_\alpha,...,y_{|\mathcal{C}|} :y_\alpha \in \mathbb{N}^{|\mathcal{S}|}\}\\
    \mathcal{R} &= \{ \ce{ $y_{\alpha}$ ->[$k_{y_{\alpha} \to y_{\beta}} \cdot x^{y_{\alpha}}$] $y_{\beta}$ } : k_{y_{\alpha} \to y_{\beta}} \geq 0\},
\end{align*}
where Roman letters ($i,j$) and Greek letters ($\alpha,\beta$) are used to
denote species and complex indices, respectively. A complex is a multi-set
of species, and is described by a column vector $y_\alpha$ representing the
stoichiometry/direction of the multi-set. The actual state of the system is
described by the vector $x=(s_1, ... , s_i,..., s_{|\mathcal{S}|})$ (that
is, a species and its quantity is denoted the same way). The term
$x^{y_\alpha}$ is defined as
\[
x^{y_\alpha}=\prod_{i=1}^{|\mathcal{S}|} {x_i}^{y_{\alpha,i}}.
\]

\eeD
  \beXa [reactions representation of SIRS  without demography]
{
\be{RNSIR} \bc  \ce{ $s$ + $i$ ->[\beta $s i$] $2i$}\\ \ce{ $i$ ->[\gamma_i $i$] $r$ } \\
 \ce{$r$ ->[\gamma_{r} $r$] $s$ } \\
 \ce{  $s$  ->[\gamma_{s} $s$] $r$ } \ec \ee
}

Here, \begin{align*}
    \mathcal{S} &= \{s,i,r\}, x=(s,i,r),\\
    \mathcal{C} & = \{y_1=(1,1,0)^t,y_2=(0,2,0)^t,y_3=(1,0,0)^t,y_4=(0,1,0)^t,y_5=(0,0,1)^t \}\\
    \mathcal{R} &= \{ \ce{ $y_1$ ->[$k_{y_1 \to y_2}x^{y_1} = \b
        x^{y_1}=\beta s i$] $y_2$},
    \ce{ $y_4$ ->[$k_{y_4 \to y_5}x^{y_4} = \gamma_ i  x^{y_4}=\gamma_ i  i$] $y_5$},...\},
\end{align*}
\eeXa

\subsection{Locking sets and semilocking sets/siphons}

Recall first that the zero coordinate set $Z(s)$ of any fixed boundary point in a mass action reaction network is a {\emph siphon/semilocking set} \cite[Prop. 4.5]{AdLS}, \cite{ShiuStu},
\cite[Lem. 63]{FreFelW}, a concept introduced in the CRN literature by Anderson \cite{AndGAS} and  Angeli et al.\ \cite{AdLS}, where the last authors made also the connection with the parallel concept of siphons in the Petri networks literature (see also \cite{koch2010modeling} for a survey of the Petri net
approach in systems biology). Siphons  provide essential information for \nne\ ODEs about whether trajectories can escape from or remain trapped on boundaries.


\beD[siphon/semilocking set, locking set]\lbl{d:sip}\cite{feinberg1987, Ang07}.
 \leavevmode

 \begin{itemize}
\item A \textbf{siphon/semilocking set} $W \subseteq \mathcal{S}$ is a nonempty subset of species such that whenever a species in $W$ appears in a product complex, at least one species in $W$ must appear in the corresponding reactant complex.

\item A \textbf{locking set} $W \subseteq \mathcal{S}$ is a nonempty subset where every reaction has at least one species from $W$ in its reactant complex.

    \im A siphon/semilocking set is \textbf{minimal}  when it contains no other siphon included within.
\end{itemize}

\eeD

\beR Note that  the RN representation is not unique, and that  the search for ``good representations" is a very active field -- see for example \cite{Buxton,Hong} (whether the ``good representations" in chemistry, which enjoy weak-reversibility and low-deficiency, for example, will also be useful in ME is a topic for further research).  However, the only theoretical CRN concept used in this paper, the set of minimal siphons, depends only on the ODE and not on the RN representation -- see Remark \ref{r:inv}.\eeR

\beR \label{r:inv}
Angeli, De Lenheer, Sontag \cite[Prop 2]{AdLS} (see also   Shiu and Sturmfels \cite[Prop 2.1]{ShiuStu}) proved that a nonempty set $I$  is semilocking/siphon iff the boundary face $F_I$ is forward-invariant for the dynamics. This implies that siphons are determined by the ODE, and are the same for any RN representation (while the proof of \cite[Prop 2]{AdLS} does use a specific representation, any other representation would have lead to the same result).
\eeR

The  semi-locking  property is  fundamental for understanding persistence.
\beD[persistence]\label{d:per} For a mass--action system  on $\mathbb{R}^n_{\ge 0}$, \emph{persistence} means trajectories starting in the positive orthant remain bounded away from the boundary: if $x_i(0)>0$ for all $i$, then $\liminf_{t\to\infty}x_i(t)>0$ for all $i$ \cite{AdLS}.
\eeD

 \cite{AdLS} showed that if the $\omega$-limit set does not intersect semi-locking boundaries except at equilibria, and if all trajectories starting on non-semilocking boundaries eventually leave those boundaries,  then the system exhibits persistence.

\ssec{Constructing  minimal siphons for ME models, by eye inspection}
ME models are characterized by the presence of  infection-type reactions, like S+I$\to$ 2 I, which are called in CRN ``autocatalytic".
\beD[autocatalytic reactions]\lbl{d:aut}\cite{Horn72}
 A reaction is called autocatalytic if there exist a species which appears
both as reactant and as product, with net production (coefficient as product exceeds that as reactant).
\eeD

\beR \BI \im  An input species (which has external arrivals) cannot be included in any siphon. \im  A species which is product of an autocatalytic reaction is always included in a (minimal) siphon.
\im The total set of species, excluding the input species,  is always a locking set (and siphon).
\EI
\eeR

The example below  introduces the reader to the CRNT idea that  crucial algebraic information like  siphons may be obtained 
 directly by examining the RN (first) column of the reaction description.

\beXa [detecting the siphons of an example with two ``infections", one ``combination" and two `` mutations"]
The   model    defined by the
reactions and rates representation below
has  two intersecting minimal siphons.

\noindent
\begin{minipage}[t]{0.35\linewidth}
\begin{verbatim}
RN = {0 -> "S",
 "S" -> 0,
 "I1" -> 0,
 "I2" -> 0,
 "I3" -> 0,
 "S" + "I1" -> 2*"I1",
 "S" + "I2" -> 2*"I2",
 "I1" + "I2" -> "I3"
 "I3" -> "I1",
 "I3" -> "I2",
 };
\end{verbatim}
\end{minipage}\hspace{0.03\linewidth}%
\begin{minipage}[t]{0.40\linewidth}
\begin{verbatim}
rts = {Lambda,         (* birth *)
 muS*S,                (* S death *)
 mu*I1,               (* I1 death *)
 mu*I2,               (* I2 death *)
 mu*I3,               (* I3 death *)
 ba1*S*I1,             (* infection I1 *)
 be2*S*I2,             (* infection I2 *)
 de*I1*I2}             (* combination *)
 m1*I3,               (* mutation I3→I1*)
 m2*I3,               (* mutation I3→I2*)
 ;
\end{verbatim}
\end{minipage}\\

with ODE
$$X'=\bep-\beta _1 i_1 S-\beta _2 i_2 S-\beta _3 i_3 S+\lambda -S \mu _s\\-\delta  i_1 i_2-i_1 \mu _1+i_3 m_1+\beta _1 i_1 S\\-\delta  i_1 i_2-i_2 \mu _2+i_3 m_2+\beta _2 i_2 S\\\delta  i_1 i_2-i_3 \mu _3-i_3 m_1-i_3 m_2+\beta _3 i_3 S\eep$$
\eeXa
\beR  ``Combination" reactions like $"I1" + "I2" \rightarrow "I3"$ are absent from classical epidemic models, but may occur in viral models (for example, a healthy cell and a virus combine into one unhealthy cell).\eeR

Here, after excluding S, which may not be involved in any siphon due to its external arrivals, we find that the  species $\{I_1\}$ and $\{I_2\}$ which intervene in autocatalytic reactions may be completed to the minimal siphons $\{I_1,I_3\}$, and $\{I_2,I_3\}$. It may be checked that the NGM matrices have the following triangular structure:
$$F=\left(
\begin{array}{ccc}
 \beta _1 S & 0 & 0 \\
 0 & \beta _2 S & 0 \\
 0 & 0 & 0 \\
\end{array}
\right), V=\left(
\begin{array}{ccc}
 \mu _1 & 0 & -m_1 \\
 0 & \mu _2 & -m_2 \\
 0 & 0 & \mu _3+m_1+m_2 \\
\end{array}
\right), K=\left(
\begin{array}{ccc}
 \frac{\beta _1 S}{\mu _1} & 0 & \frac{\beta _1 m_1 S}{\mu _1 \left(\mu _3+m_1+m_2\right)} \\
 0 & \frac{\beta _2 S}{\mu _2} & \frac{\beta _2 m_2 S}{\mu _2 \left(\mu _3+m_1+m_2\right)} \\
 0 & 0 & 0 \\
\end{array}
\right)$$

We note here the usual upper triangular structure, which implies that
$R_0=\max[\frac{\beta _1 S}{\mu _1},\frac{\beta _2 S}{\mu _2}],$ in a model which does not have a straightforward ME interpretation.
\medskip

\subsection{From critical minimal siphons to the disease free equilibrium\ (DFE), and locally monotone x-y models}

The only minimal siphons  encountered in actual ME models have the property that their supporting boundary faces  can be either attracting or repelling (depending on
the parameters, or on reproduction and invasion numbers/rates). Cf. \cite{AdLS}, such siphons must be  ``critical" in the sense of satisfying the following:
\beD[critical siphon]\cite{AdLS}.
A siphon $T$ in a CRN with \sm\ $\Gamma$ is \emph{critical} if it contains no
support of a positive conservation relation,  that is, if
 there exists no nonzero vector $c\ge0$ with
$c^\top\Gamma=0$ and $\operatorname{supp}(c)\subseteq T$.
\eeD

\beR A) The concept of critical siphon is also dependent only on the ODE, and not on a particular
RN representation.

B) Boundary faces associated with uncritical siphons may not be reached by the ODE, when starting in their exterior \cite{AdLS}. \eeR

Since minimal siphons are considerably easier to compute
than \fbp s (and the chemical reaction networks literature provides several algorithms for
that),  they are the first thing we compute  for any
autocatalytic model. For small models, this may be achieved by eye inspection,
as in the example below.

\beXa[Example with  both critical and non–critical minimal siphons]
\leavevmode

\noindent
\begin{minipage}[t]{0.43\linewidth}
\begin{verbatim}
RN = {
  0 -> "S",
  "S" -> 0,

  "S" + "I1" -> 2*"I1",
  "S" + "I2" -> 2*"I2",

  "D" -> "E",
  "E" -> "D"
};
\end{verbatim}
\end{minipage}\hspace{0.03\linewidth}%

There are three minimal siphons:
\(
T_1=\{I_1\}\ (\text{critical}),\quad
T_2=\{I_2\}\ (\text{critical}),\quad
T_3=\{D,E\}\ (\text{non--critical}).
\)
 \(T_1\) and \(T_2\), induced by the \aut\ reactions, include no positive conservation, hence they are \emph{critical}.
The siphon  with reversible complex graph \(\{D,E\}\) supports the positive linear conservation \(D+E=\text{const}\), hence it is \emph{non–critical}.
\eeXa

The minimal critical siphons (whose associated invariant faces   may  be reached by the ODE, and have maximal dimension), provide a rigorous definition of the disease free equilibrium (DFE).

\beD [DFE] \label{d:DFE} The union $TS=\cup_E I_{E}$ of the indices  of all {\bf minimal critical siphons} $E$ will be called the total siphon, and the face $W_{TS}=\{X:X_i=0, \forall i \in TS\}$ will be called the DFE face.

If a boundary fixed point $E_0\in W_{TS}$ exists,  it  will be called DFE.
\eeD
\beR The total siphon is itself an invariant face, by Remark \ref{r:inv}. \eeR

\beR
\leavevmode
\noindent

\BI \im We might  restrict the definition \ref{d:DFE} to all minimal critical siphons which are both
self-replicable and drainable (see \cite{DesGop} for definitions),  since minimal siphons are always both
self-replicable and drainable in ME examples.
However, we do not do this now, for simplicity.

\im Minimal self-replicable siphons are  connected conceptually to the   concept of cores/minimal self-amplifying subhypergraphs \cite{BlaGag}, which originate in the exciting theory of auto-catalytic systems -- see  \cite{BlokLacNghe,PengGag,VasSta,GagPoly}, and  see section \ref{s:rel} for some further definitions. An alternative definition for the DFE zero set could be as the union of all the cores (in our examples this coincides with the union of the minimal siphons).

Note that precise relations between the two concepts do not seem to have published yet; however, they suggest  relations between multi-strain ME models and the
 multi-tier Seed-Dependent Autocatalytic Systems (SDAS) \cite{Peng2021SDAS,Peng2022HierarchicalSDAS,Peng2023ChemicalEcosystem,
 Blok23} --see section \ref{s:SDAS}. Relating the two fields  seems a  promising research direction.

\im If there exist several boundary fixed points whose  set of zero values is $TS$, they might be all called DFEs, but for now we are excluding this possibility in our definition \ref{d:ME} of ME models.
\EI
\eeR

\beD  [x-y models, Lotka-Volterra-Kolmogorov/generalized Lotka-Volterra models] \label{d:xyLV} $~$

Here, at  each boundary fixed point, the ``invading  variables" $\x$ are those that take
the value $0$, and the ``resident  variables" $\y$ are their complement.

A) A positive ODE which admits  a DFE $E_0$, will be called an x-y model. In this case, after  removing eventual conservations,
  the ODE \mbw:
\begin{equation}
    \label{de2}
    \left\{ \; \; \;
    \begin{aligned}
    \frac{d\x}{dt} & = \mathbf{f}(\x,\y)= M(\x,\y) \x, \\
    \frac{d\y}{dt} & = \mathbf{g}(\x,\y).
    \end{aligned}
    \right.
\end{equation}
where $\x,\y$ denote the zero and positive coordinates of the DFE.

B)  x-y models with   equations  of the form
\begin{equation}\label{eq:LV}
  x_i'  = x_i f_i (\x,\y), \for i\in \x_{E_0},
\end{equation}
\ie\ with $M$ diagonal,   will be called below \LV-Kolmogorov-type (LVK)
x-y models.

When   $\y$ further varies in a compact space, these  were
studied  by \cite{PatelSchreiber,HofSch}.

C) An $\xy $ model with $\y =\emptyset$ will be called $\x$ model.

D) An $\x $ model with Kolmogorov-type equations
for each $\x$ \var\
\begin{equation}\label{eq:cLV}
  x_i'  = x_i f_i (\x), \for i\in \x_{E_0},
\end{equation}
will be called \gLV\ (GLV) following Goh and subsequent literature, or \LVK\ (LVK).

\eeD

\beR  The matrix $M(\x,\y)$ \sats, by \eqr{de2},   that the non-zero steady states  must all satisfy
 \be{M} det(M(\x,\y))=0,\ee
 but is not uniquely defined.

 For typical ME models, it is possible to choose a  matrix $M(\x,\y)$ appearing in the   equations for $\x'$ at the DFE which is  Metzler, for any $\y \in \mathbb{R}_{\geq 0}^{m}$  which is  equivalent to saying that the flow of the $\x$ \var s, when the other \var s are fixed,  is monotone.

\eeR

One foundational result of ME, the NGM method \cNGM\ reviewed in
Section~\ref{s:NGM}, consists in expressing the stability domain
of the disease free equilibrium\ in the form $R_0\leq 1$. This result is related to  the regular splitting theorem reviewed in next section.

\subsection{Defining ME models by assuming the existence of a regular splitting for the Jacobian $J_x$ and stability of the Jacobian $J_y$ at the DFE} \label{sec:heuristicNGM}

Currently, there is no accepted definition of ME models. It is accepted  that ME models are a subset of mass-action representable \nne\ models, but this is a very general class. We believe  that ``relevant"  ME models are  a subset of CRN models, a class which in itself rules out certain ``physically unrealistic behaviors like hidden inflow reactions (HIR)/cornucopias \cite{MulFlaSta} like
$D+Y \rightarrow 2 D + Y $ (which result into a \nne\ stoichiometric vector, but are not pure inflow reactions).

\Fr we conjecture that the existence of a regular splitting  for the $x-y$ partition induced by the DFE could be taken as definition of ME models, and we offer two possible definitions in this direction.

\beD [ME model,    classic and algorithmic] \label{d:ME}
\noindent
\BI \im We will call classic ME model  a CRN with  no hidden inflow reactions, and  such that:

  \BEN \im The DFE index set (the union  all minimal critical siphons) is nonvoid, and contains a compact invariant set, which is composed of precisely one point;  furthermore, this point
   may be, depending on the \para s, either stable or unstable.
    \im   The Jacobian $J_y$ at the DFE, where $\y$ are the non-DFE (resident) \var s,  is stable.
    
    \im  The Jacobian $J_x$ at the DFE is Metzler (equivalently, the ODE is locally  increasing in $\y$ at the DFE)  \im  $J_x$ admits some regular splitting  $J_x=F-V$ (and as a consequence, the instability domain of the DFE may  be expressed as $R_0=\rho(F V^{-1})>1$).

    \EEN
\im We will call algorithmic ME model an  x-y model satisfying 1. - 4. above, and such that 5. is \satd\ with F  obtained from the totality of reactions involving x and y reactants, by removal of all terms preceded by a minus, as proposed in \cite{AABJ}.
\EI
\eeD

\sec{A two strain SI$^2$V Lotka-Volterra-Kolmogorov model with permanent cross-immunity and vaccination, where the persistence conjecture  holds \cite{RahVac}; SI2V.nb}\label{s:Rah}
This section presents an example which \sats\  both the persistence conjecture and the \CEP.
\beD[Competitive Exclusion/Linear Complementarity Principle]
\label{d:LCP}
We will say that a k-strain ODE \sats\ the \CEP\ (CEP) or Linear Complementarity Principle   if its parameter space
may be partitioned the  into $2^k$ stability regions, each of which  contains precisely one LAS equilibrium.
\eeD

The SI$^2$V two-strain  model  of \cite{RahVac} is defined by the following
reactions + rates representation:
\be{RNv}\left(
\begin{array}{cc}
 0\to  {s} & \Lambda(s)  \\
  {i_1}+ {s}\to 2  {i_1} & \beta _1 i_1 s \\
  {i_2}+ {s}\to 2  {i_2} & \beta _2 i_2 s \\
  {s}\to  {v} & \rho  s \\
  {i_2}+ {v}\to 2  {i_2} & i_2 v \beta _v \\
  {v}\to 0 & \mu_v   v \\
  {i_1}\to 0 & i_1 \mu_1 \\
  {i_2}\to 0 & i_2 \mu_2 \\
\end{array}
\right)
\ee
(this is the simplified version, where the recovered class, which does not affect the others, is omitted).

For a discussion of the advantages of using this representation for ME models, see \cite{VAA,AAN,AAHK}.

The two minimal siphons are {{$i_1$},{$i_2$}}.
The induced ODE is:
\be{ODERahv} X'=\bc-\beta _1 i_1 s-\beta _2 i_2 s+\Lambda( s)-\rho  s\\
\bep \mR _1  s- 1&0\\0&\mR _2  s+  \mR _v v- 1\eep
\bep \mu_1 i_1 \\\mu_2 i_2\eep\\
-i_2 v \beta _v+\rho  s-\mu_v  v\ec:=\bc-\beta _1 i_1 s-\beta _2 i_2 s+\Lambda( s)-\rho  s\\
M(  s, v)
\bep \mu_1 i_1 \\\mu_2 i_2\eep\\
-i_2 v \beta _v+\rho  s-\mu_v  v\ec
\ee
with
\[
{\mR}_j=\frac{\beta_j}{\mu_j}, j=1,2,\qu
{\mR}_v=\frac{\beta_v}{\mu_2}.\]

\beA From now on, we will assume that $\La(s) =\La - \mu s:=\Lambda(1 -s /s_0)$.
  \eeA

 \beR

The fact that in \eqr{ODERahv} the $M=Diag[f_j(\y)]=Diag[R_j(\y)-1], \y=(s,v)$  is diagonal makes the infection problem (with fixed $s$) a Lotka-Volterra-Kolmogorov ODE (see definition \ref{d:xyLV}), with $$f_1 = \mathcal{R}_1 s - 1, f_2 = \mathcal{R}_2 s + \mathcal{R}_v v - 1.$$

 We expect therefore to find  four possible  fixed points (see also \cite{Rahman,AAHK} for a  similar problem); furthermore,
 the diagonal structure of $M$ suggests combining the local stability problems for all the fixed points  into the single  \LCP\ type  optimization problem
 defined in \eqr{LCP}.
  \eeR

  \beD[The linear complementarity problem associated to a  Lotka-Volterra-Kolmogorov ODE]
 The linear complementarity problem associated to a  Lotka-Volterra-Kolmogorov ODE is:
 \be{LCP}\bc
 max \sum_j \mu_j i_j f_j, \\
 i_j\geq 0, f_j\leq 0, \for j \ec.\ee
 \eeD

 Clearly, \eqr{LCP} obtains the maximum value $0$, and the solutions  are into a one to one correspondence with   cases
  of the type $i_S=0, f_{S^c}\leq 0,$  (where S runs over all   subsets of the infection set), which may be seen to correspond exactly to the stability conditions for  individual fixed points.

We conjecture that a unified proof for CEP/LCP, which has only been established in particular cases,  may be provided by rigorously defining k-strain models, and using their \repF, defined below:

\beD [simple k-strain model] \label{d:mtS} A simple k-strain model is an ME model
 for which the  NGM matrix  has, after permutation,  a triangular block structure with {\bf precisely $k$ blocks}, each of which has a unique unconditionally positive eigenvalue, and for which in all siphon faces which contain a compact invariant set in their interior, this is composed of precisely one point.
\eeD

\begin{corollary}
For a simple k-strain model
$\rho(K)=\max_{j=1}^k \rho(K_{jj})$, where $ k$ is the number  of  blocks
in the decomposition of the NGM.
\end{corollary}

\beD[reproduction functions, invasion numbers  and reproduction  numbers for simple k-strain models] \label{d:reF}
\leavevmode

For a simple k-strain model:
\noindent
\BEN \im the unique strictly positive eigenvalues of the NGM blocks, $R_i(\y), i=1, ...k$, with the $\y$ \var s left free, will be called {\bf R-reproduction functions}
(associated to the  strain $i$).

\im  
\be{Ri} R_i(s_0)={\mR} _i s _0 :=R_i^{\T 0}:=R_i, i=1,...,k\ee
will be called the  {\bf \brn\ of strain} $i$.

\im   Finally, in the case $k=2$
with two strains  which correspond to two minimal siphons with resident values
$s_j$,
\be{invN}R_i^{\T j}:=R_i(s_j)\ee will be called
the {\bf invasion number of invading strain $i$ on resident strain $j$} (see Definition \ref{d:inv3}  for the general case).

\EEN
\eeD

\beR Note that here we have given a precise mathematical definition
for the ME  concept of invasion number, as the reproduction function of one (invading) strain, evaluated at the resident boundary fixed point\ of the other strain (this is  formalized  to more than two strains  below). Note also that in ecology, the invasion rates turn out to be proportional to $R_i(s_j)-1$, in the simple case when resident invariant sets are singletons.
\eeR

\beXa [reproduction functions and invasion numbers of SI2V] For the \cite{RahVac} model,
the NGM $$K=\left(
\begin{array}{cc}
 \frac{\beta _1 s}{\mu_1 } & 0 \\
 0 & \frac{\beta _2 s+\beta _v  v}{\mu_2 } \\
\end{array}
\right):=\left(
\begin{array}{cc}
 \mR_1 s & 0 \\
 0 & \mR _2 s+\mR _v  v \\
\end{array}
\right) $$
 reveals that the reproduction functions --see Definition \ref{d:reF}--
 are \be {rf} R_1(s,v):=\mR _1 s, R_2(s,v):=\mR _2 s + \mR_v  v.\ee

  The partition or the \para\ space in terms of ${R}_1, {R}_2, R_2^{\T 1} ,R_1^{\T 2} $ turns out to be:

\begin{enumerate}
\item  {\bf DFE stable}: $\max[R_1, R_2] \leq 1$.
\item  {\bf $E_1$ stable}: $R_1 > 1, R_2^{\T 1} \leq 1$, i.e. strain 2 cannot invade.
\item  {\bf $E_2$ stable}: ${R}_2 > 1, R_1^{\T 2}  \leq 1$, i.e. strain 1 cannot invade.
\item  {\bf $E^*$ stable}: $\min[R_1, R_2,R_2^{\T 1} , R_1^{\T 2} ] > 1$, i.e. both  strains 1,2 can exist and can invade.
\end{enumerate}
\eeXa
\beR  DFE is stable iff it is the only boundary equilibrium.
The coexistence equilibrium is stable iff both single-strain equilibria   exist and are unstable.\eeR

  \subsection{The  Fixed Points, and the Solution of the Linear Complementarity Problem}

The DFE, E1, and the coexistence point EE turn out to be explicit, and the solution of E2 reduces to resolving a quadratic. Finally, the LCP is equivalent to the following four cases, obtained by letting $S$ run over the two minimal siphons, their union, and the empty set:

(1) $S=\{i_1,i_2\} \Eq i_1 = i_2 = 0$ (with $f_1, f_2 \leq 0$ for stability) yields the
{\bf Disease-Free Equilibrium (DFE)} $E_0$:
$$E_0 = \left(s_0 = \frac{\Lambda }{\mu +\rho }, \, 0, \, 0, \, v_{0} = \frac{\rho }{\mu_v} s_0\right),$$
and  direct or  NGM  analysis confirm the expected result that stability holds iff $R_0 =\max [R_1,R_2]\leq 1,$
where $R_1,R_2$ are the reproduction numbers of each strain --see definition \ref{d:reF}.3.

(2) $S=\{i_2\} \Eq i_1 > 0, i_2 = 0 \Rightarrow f_1 = 0, f_2 \leq 0$  yield the
    {\bf strain 1 only} equilibrium:$$E_1 = \left(s_1 = \frac{\mu _1}{\beta _1}, \, i_{1,1} =\frac{\beta _1 \Lambda -\mu _1 \rho -\mu  \mu _1}{\beta _1 \mu _1}= \frac{\rho+\mu}{\beta_1}(R_1 (s_0)- 1)=\frac{\rho+\mu}{\beta_1 s_1}
     (s_0 - s_1), \, 0, \, v_1 = \frac{\rho }{\mu_v} s_1\right),$$
with existence condition:
\be{R1}R_1 :=R_1 (s_0) > 1 \Eq s_0 > s_1,\ee
where we recall that the R-reproduction function $R_1 (s)$ is defined in definition \ref{d:reF}.1, as the unique strictly positive eigenvalue of the NGM, with the y variables left free, associated to the block corresponding to strain $1$.

The analysis of the Jacobian at $E_1$ (see Mathematica file)  confirms the predicted stability condition
 $$\Rot:=R_2(s_1,v)=\mathcal{R}_2 s_1 + \mathcal{R}_v v \leq 1.$$

(3) $S=\{i_1\} \Eq i_2 > 0, i_1 = 0$, the {\bf strain 2 only equilibrium} has $i_2$ \satg\ a quadratic equation $A i_2^2+Bi_2+C=0$, with \coe s:
$$ C=\mu _2 (\mu +\rho ) \mu _v -\Lambda(\left(\beta _2   \mu_v +   \rho  \beta _v\right)), B=\mu _2 (\mu +\rho ) \beta _v+
\beta _2 \left(\mu _2 \mu _v-\Lambda  \beta _v\right), A=\beta _2 \beta _v \mu _2 $$

As proved in \cite{RahVac}, A positive free coefficient $C$ occurs
iff \be{R2}R_2 :=R_2 (s_0,v_0) > 1,\ee
and in this case the quadratic has precisely one positive root.
\Fr $s,v$ are also positive if $i_2$ is.

The direct LAS approach involves analyzing the stability of a third order polynomial, and  reveals that the result is as expected, stability holds iff $\Rto<1$.

(4) $S=\emptyset$ yields the
{\bf endemic equilibrium $E^*$}, which,  surprisingly, is rational

$$s= s_1=\frac{\mu _1}{\beta _1},
v= \frac{\mu _1 \mu _2 (\mR _1 -\mR _2)}{\beta _1 \beta _v},
i_2= \frac{ \rho }{\mu _2(\mR _1 -\mR _2)}-\frac{\mu _v}{\beta _v},
i_1= \frac{\Lambda }{\mu _1}+\frac{\frac{\beta _2 \mu _v}{\beta _v}-\mu }{\beta _1}-\frac{\rho }{\mu _1(\mR _1 -\mR _2)}.$$

It turns out  that the (positivity) existence \con s
may be written in the form
\be{exL}\mR_1 >\mR_2>1, \Rot >1,\Rto >1;\ee
the first three are obvious, and the last two may be checked easily to be  the positivity conditions for $i_1,i_2$ (or, see file  SI2V.nb).

\eqref{exL}  are precisely the expected LAS conditions.
That this is indeed the case, namely that existence implies LAS could be proved by the direct LAS approach, but this    is rather challenging, since it requires solving 4'th order \RH\ conditions. Fortunately, LAS has been already established in \cite{RahVac}, by offering a \Lf\ which shows that GAS holds whenever EE exists.
\beT[persistence of the SI2V model]  The SI2V model of \cite{RahVac} satisfies both the persistence conjecture and the \CEP.
\eeT

We conclude this section with  an open problem.

\beO [LCP/CEP for k-strain models] \lbl{o:LCP} Which k-strain models \saty\ LCP/CEP, and under what conditions the   stability regions for the fixed points may be  expressed in terms of the $2^k$ reproduction functions $R_{\T S^c}^S$, where $S$ runs over the $k$ minimal siphons?
\eeO

\section{A directed graphs on minimal siphons (IGMS)}\lbl{s:IGMS}
This section establishes that  next–generation matrices \(K=F V^{-1}\) have a block–triangular structure from $V$ whenever the minimal siphon partition the infected, and the graph defined below is acyclic.

\beD[The Interaction Directed Graph on Minimal Siphons (IGMS)]\lbl{d:IGMS} The Interaction Directed Graph on Minimal Siphons (IGMS) is the graph on the minimal siphons $T_j, j=1, ...,m$, with a directed edge $T_i \to T_j$ whenever there exists a reaction net-producing at least one species in $T_j$ from reactants in $T_i$.
\eeD
\beXa  [a  tree IGMS for a three-tier Seed-Dependent
Autocatalytic Systems (SDAS) -type  ME model]
The RN
$$\left(
\begin{array}{cc}
 0\to \text{S} & \lambda  \\
 \text{S}\to 0 & S \mu _s \\
 \text{I1}+\text{S}\to 2 *\text{I1} & \beta _1 i_1 S \\
 \text{I1}+\text{I2}\to 2 * \text{I2} & \beta _2 i_1 i_2 \\
 \text{I2}+\text{I3}\to 2 *\text{I3} & \beta _3 i_2 i_3 \\
 \text{I1}\to 0 & i_1 \mu _1 \\
 \text{I2}\to 0 & i_2 \mu _2 \\
 \text{I3}\to 0 & i_3 \mu_3 \\
\end{array}
\right)$$
has three disjoint minimal siphons, and a tree IGMS, as shown in Figure \ref{fig:strain_tree}. The NGM matrix is triangular (in fact diagonal), in line with theorem \ref{t:par}.
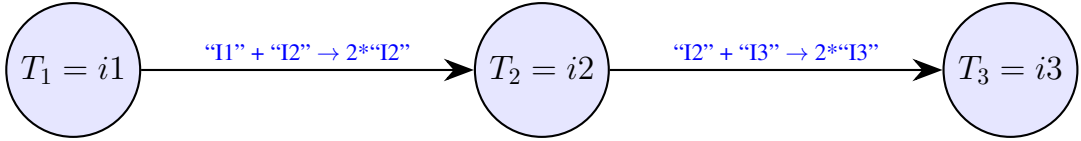
\begin{figure}[h!]
\centering
\begin{tikzpicture}[
  scale=0.89,  
  transform shape,  
  strain/.style={draw, circle, minimum size=2cm, thick, fill=blue!10, font=\Large},
  arrow/.style={-{Stealth[length=4mm]}, thick, >=stealth},
  reaction/.style={font=\small, color=blue}
]

\node[strain] (T1) at (0,0) {$T_1 = i1$};
\node[strain] (T2) at (7,0) {$T_2 = i2$};
\node[strain] (T3) at (14,0) {$T_3 = i3$};

\draw[arrow] (T1) -- node[above, reaction, pos=0.5] {``I1'' + ``I2'' $\to$ 2*``I2''} (T2);
\draw[arrow] (T2) -- node[above, reaction, pos=0.5] {``I2'' + ``I3'' $\to$ 2*``I3''} (T3);

\end{tikzpicture}
\caption{Three-tier  sequential SDAS-type  ME model. Siphons $T_i$ are represented
by circles. The reactions are of  autocatalytic infection type, where contact between strains leads to conversion.}
\label{fig:strain_tree}
\end{figure}

This RN can not be interpreted as an existing ME  model, but it has been used in the quite similar Seed-Dependent Autocatalytic System (SDAS) theory. Originally proposed to model stepwise prebiotic chemical evolution, this is now used to identify hierarchically nested autocatalytic motifs in large reaction networks \cite{Peng2022HierarchicalSDAS,Peng2023ChemicalEcosystem,BlokLacNghe,
Blok23}. This theory provides general results for ecological  and potential multi-strain epidemic systems, where new functional tiers emerge only after the appearance of specific seed species.
\eeXa

\subsection{Block structure of the Next--Generation Matrix under Acyclic Minimal Siphon Decompositions (AMSD)}

We show in this section
that in epidemic models where the minimal siphons partition the DFE support and IGMS  is acyclic, we have: (i) \emph{block–triangular} NGMs $K=FV^{-1}$, after reordering,  and (ii) \emph{decomposition of the threshold} $\rho(K)=\max_i \rho(K_{ii})$ by blocks.

\beT[acyclic minimal siphon decomposition implies triangularity of NGM]\label{t:par}
Let $\mathcal{I}$ be the set of infected variables and suppose $\{T_1,\dots,T_m\}$ is a decomposition of $\mathcal{I}$ into siphons, and let $G_{\mathrm{out}}$ denote the IGMS  (on vertices $\{T_1,\dots,T_m\}$, with  a directed edge $T_i \to T_j$, $i\neq j$, whenever there exists a reaction whose reactant set intersects $T_i$ and whose product set intersects $T_j$).
Assume $G_{\mathrm{out}}$ is acyclic.
Let $\pi$ be any topological ordering of $G_{\mathrm{out}}$, and reorder the infected variables by listing first all species in $T_{\pi(1)}$, then those in $T_{\pi(2)}$, and so on, with arbitrary order within each block.
Under this ordering, $F$ is block diagonal and $V$ is block lower triangular; hence $K = F V^{-1}$ is block lower triangular with $m$ diagonal blocks.
If $G_{\mathrm{out}}$ has no edges, then $K$ is block diagonal.
\eeT

\begin{proof}
Let $\mathcal{I} = T_1 \sqcup \cdots \sqcup T_m$ be the decomposition of infected variables into minimal siphons, and let $\pi$ be a topological ordering of $G_{\mathrm{out}}$. Write the linearized infection subsystem as $\dot x = (F - V)x$,
where $F$ represents new infection terms and $V$ the transition and removal terms.

The acyclicity of the IGMS implies that both  $F$ and  $V$ are block lower triangular, and hence
 the product $K = FV^{-1}$ is also block lower triangular.

\end{proof}

\begin{corollary}
If $\mathcal{I}$ admits an acyclic siphon decomposition, then $\rho(K)=\max_i \rho(K_{ii})$ after a compatible ordering of variables; hence threshold behavior decomposes by blocks.
\end{corollary}

 \beR AMSD is a sufficient, but not necessary condition for block triangular NGM structure, as the next example shows.
 \eeR

\ssec{Example of cyclic IGMS, with block diagonal NGM structure: 5cycles.nb}
\noindent
\begin{minipage}[t]{0.46\linewidth}
\begin{verbatim}
RN = {
  0 -> "S",
  "S" -> 0,
  "S" + "I1" -> 2*"I1",
  "S" + "I2" -> 2*"I2",
  "I1" + "I2" -> "I12",
  "I1" + "I12" -> "I2",
  "I2" + "I12" -> "I1",
  "I1" -> 0,
  "I2" -> 0,
  "I12" -> 0
};
\end{verbatim}
\end{minipage}\hspace{0.03\linewidth}%
\begin{minipage}[t]{0.48\linewidth}
\begin{verbatim}
rts = {
  la,             (* birth *)
  mu*S,           (* S death *)
  be1*S*I1,       (* inf. I1 *)
  be2*S*I2,       (* inf. I2 *)
  de*I1*I2,       (* I1+I2->I12 *)
  et1*I1*I12,     (* I1+I12->I2 *)
  et2*I2*I12,     (* I2+I12->I1 *)
  mu1*I1,         (* I1 death *)
  mu2*I2,         (* I2 death *)
  mu12*I12        (* I12 death *)
};
\end{verbatim}
\end{minipage}

\noindent
The minimal siphons: \(T_1=\{I_1,I_{2}\}, T_2=\{I_2,I_{12}\}, T_3=\{I_1,I_{12}\}\) yield
an IGMS with two three–cycles \(T_1\rightleftarrows T_2 \rightleftarrows T_3\),
and \(T_1\rightleftarrows T_3 \rightleftarrows T_2\), and three pair cycles. Despite this, the NGM has block-diagonal structure:
$$K=\left(
\begin{array}{ccc}
 \frac{\beta _1 s}{\mu _1} & 0 & 0 \\
 0 & 0 & 0 \\
 0 & 0 & \frac{\beta _2 s}{\mu _2} \\
\end{array}
\right)$$

\section{A  two strain model with temporary cross-immunity, with ADE and immunity waning,    that may exhibit Hopf bifurcations:~\cite{ferguson1999effect,schwartz2005chaotic,
nuno2005dynamics,Agu07,Chung,GavishMusa,GavCEP,Gavnew};GavScan.nb}\label{s:Gav}
\ssec{Background}
The two strain model may be traced back to \cite{ferguson1999effect,schwartz2005chaotic}, and
 is appropriate for modelling simultaneous epidemics with different pathogens, like for example Dengue and Zika. Subsequently,   two-strain models which add  further \com s allowing for temporary cross-immunity have been developed in the works of Aguiar, Stollenwerk and Kooi \cite{Agu07,Agu08,aguiar2009torus,
stollenwerk2017hopf,aguiar2022mathematical},
   and \cite{billings2008vaccinations,Bulh} examined the effects of single-strain vaccination on the dynamics of an epidemic multi-strain Dengue model (see also
\cite{AABJ} for a first public notebook). This model has also been used for several 
strains of pathogens (without immunity-effectors compartments), and in ecology 
\cite{minayev2009improving,lazebnik2023computational}.

\ssec{The \RN\ representation}

The goal of this section is both to present an interesting variant model, and also to provide some  comments on  its analysis using EpidCRN.

Figure~\ref{fig:Gav} shows a model with eight compartments: susceptibles
(s), those infected first by strain $j$ ($i_j$, primary infection), those
recovered from strain $j$ ($r_j$, as a result of primary infection), those
infected with strain $j$ after they had recovered from strain~$ i$
($y_i$) and
those recovered from both strains ($r_{12}$). There is one inflow into $S$, 8
outflows out of all compartments, and 15 inner reactions for the ``closed
model", with RN:

\vspace{-\baselineskip}
\begin{center}
\resizebox{\textwidth}{!}{\begin{tabular}{c|c|c}
Reactions & Type & Parameters \\
\hline
\makecell{
\texttt{RNc=\{ "S"+"i1"->2"i1", "S"+"i2"->2"i2",} \\
\texttt{"R2"+"i1"->"i1"+"i21","R1"+"i2"->"i2"+"i12",} \\
\texttt{"S"+"i21"->"i21"+"i1","S"+"i12"->"i12"+"i2"} \\
\texttt{"R2"+"i21"->2"i21", "R1"+"i12"->2"i12"} \\
\texttt{"i1"->"R1","i2"->"R2","i21"->"R12","i12"->"R12",} \\
\texttt{"R1"->"S", "R2"->"S", "R12"->"S" \}}
} &
\makecell{primary infection \\ cross infection \\ catalytic infection \\
 amplification \\ recovery \\ waning immunity}
&
\makecell{$\beta_1,\beta_2$ \\
$\beta_1\sigma_1,\beta_2\sigma_2$ \\
$\beta_1\eta_1,\beta_2\eta_2$ \\
$\beta_1\sigma_1\eta_1,\beta_2\sigma_2\eta_2$ \\
$\gamma_1,\gamma_2,\gamma_1,\gamma_2$ \\
$\theta_1,\theta_2,\theta_{12}$}
\end{tabular}}
\end{center}
\vspace{-\baselineskip}

\begin{figure}[H]
    \centering
\begin{tikzpicture}[
    block/.style={rectangle, draw, minimum width=1cm, minimum height=1cm, thick},
    redblock/.style={block, fill=red!80},
    blueblock/.style={block, fill=blue!60},
    gradientblock/.style={block, left color=red!80, right color=blue!60},
    arrow/.style={->, >=Stealth, thick, blue!70},
    cross/.style={->, >=Stealth, dashed, thick, blue!70},
    recovery/.style={->, >=Stealth, thick, magenta},
    catalytic/.style={->, >=Stealth, thick, green},
    amplification/.style={->, >=Stealth, thick, orange!90},
    dashedbox/.style={dashed, thick, black}
]
\node[block, fill=gray!20] (S) at (0,0) {\Large S};

\node[redblock, above right=0.8cm and 2cm of S] (i2) {\Large $\mathrm{i}_2$};
\node[redblock, right=1.5cm of i2] (R2) {\Large $\mathrm{R}_2$};
\node[gradientblock, right=1.5cm of R2] (i21) {\Large $\mathrm{i}_{21}$};

\node[blueblock, below right=0.8cm and 2cm of S] (i1) {\Large $\mathrm{i}_1$};
\node[blueblock, right=1.5cm of i1] (R1) {\Large $\mathrm{R}_1$};
\node[gradientblock, right=1.5cm of R1] (i12) {\Large $\mathrm{i}_{12}$};

\node[gradientblock, right=2cm of $(i21)!0.5!(i12)$] (R12) {\Large $\mathrm{R}_{12}$};

\draw[arrow] (S) -- node[below left, xshift=0.4cm] {$\beta_2$} (i2);
\draw[arrow] (S) -- node[above left, xshift=0.4cm] {$\beta_1$} (i1);

\draw[recovery] (i2) to[out=0, in=180] node[above, yshift=0.1cm] {$\gamma_2$} (R2);
\draw[recovery] (i1) to[out=0, in=180] node[below, yshift=-0.1cm] {$\gamma_1$} (R1);

\draw[cross] (R2) -- node[above] {$\beta_1\sigma_{1}$} (i21);
\draw[cross] (R1) -- node[below] {$\beta_2\sigma_{2}$} (i12);

\draw[recovery] (i21) -- node[above right, yshift=0.1cm] {$\gamma_1$} (R12);
\draw[recovery] (i12) -- node[below right, yshift=-0.1cm] {$\gamma_2$} (R12);

\draw[catalytic] (S) to[out=25, in=115] node[below right, xshift=0.5cm, yshift=-0.35cm] {$\beta_1\eta_1$} (i1);
\draw[catalytic] (S) to[out=-25, in=-115] node[above right, xshift=0.5cm, yshift=0.35cm] {$\beta_2\eta_2$} (i2);

\draw[amplification] (R2) to[out=0, in=90, bend right=30] node[below, xshift=0.cm] {$\beta_1\sigma_1\eta_1$} (i21);
\draw[amplification] (R1) to[out=0, in=-90, bend left=30] node[above, xshift=0.cm] {$\beta_2\sigma_2\eta_2$} (i12);

\draw[arrow, dashed, red!80] (R2) to[out=-100, in=10] node[above right] {$\theta_2$} (S);
\draw[arrow, dashed, red!80] (R1) to[out=110, in=-10] node[below right] {$\theta_1$} (S);
\draw[arrow, dashed, red!80] (R12) to[out=180, in=0] node[above right] {$\theta_{12}$} (S);

\node[anchor=north west] at (-3.5, 4) {
    \begin{tabular}{ll}
    \textcolor{blue!70}{——} & Primary infection \\
    \textcolor{blue!70}{- - - -} & Cross infection \\
    \textcolor{green}{——} & Catalytic infection \\
    \textcolor{orange!90}{——} & Amplification \\
    \textcolor{magenta}{——} & Recovery \\
    \textcolor{red!80}{- - - -} & Waning immunity
    \end{tabular}
};

\end{tikzpicture}
\caption{Schematic diagram of disease dynamics for two
    co-circulating strains. The diagram includes all reactions:
    primary infections ($\beta_i$),
    cross-infections ($\beta_j\sigma_j$),
    catalytic secondary infections ($\beta_i\eta_i$),
    amplification reactions ($\beta_j\sigma_j\eta_j$),
    recoveries ($\gamma_i$),
    and waning immunity ($\theta_i$).
    Green arrows represent catalytic
    infections where secondary cases catalyze new primary
    infections. Orange arrows represent amplification where recovered
    individuals become secondary cases through contact with existing
    secondary cases.}
\label{fig:Gav}
\end{figure}
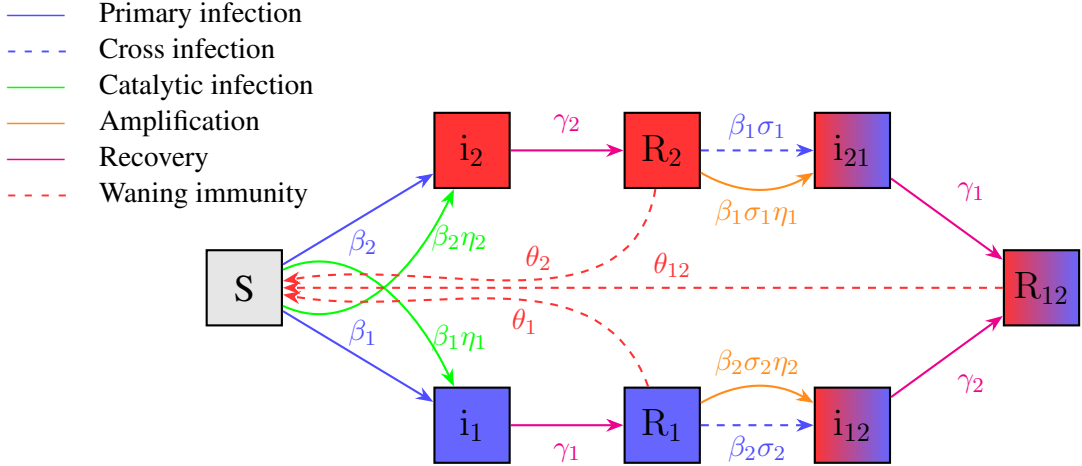

\ssec{The ODE  formulation of the  model}

We present first an extension of the mass-action model  allowing for saturation, which we plan to study in the future.
Putting \be{JF} J_k=i_k+\eta_k i_{\T k,k},
k=1,2, \T k:=3-k.\ee 
and $X=(s,i_1,i_2,i_{21},i_{12},r_1,r_2,r_{12})$, the model with saturation is given
by

\begin{equation} \label{eq:model}X'= \bc
       \f(s)-\sum_{k=1}^{2}\beta_k J_k S+\sum_{k=1}^{3}\th_k  r_k, \qu \f(s)=\Lambda-\mu s\\[1ex]
          \beta_1 J_1 s-(\mu+\gamma_1)i_1,\\[1ex]
         \beta_2 J_2  s-(\mu+\gamma_2)i_2,\\[1ex]
          \beta_1\si_{1} J_1 r_2-(\mu+\gamma_1)i_{21},\\[1ex]
          \beta_2\si_{2} J_2  r_1-(\mu+\gamma_2)i_{12},\\[1ex]
         \gamma_1i_1-\beta_2\si_{2} J_2   r_1-(\mu+\th_1)r_1,\\[1ex]
          \gamma_2i_2-\beta_1\si_{1} J_1  r_2-(\mu+\th_2)r_2,\\[1ex]
        \sum_i \ga_i y_i -(\mu +\th_{12}) r_{12},
        \ec
\end{equation}
were~$\Lambda$ is the rate at which individuals are born, $\mu 
$ is the mortality rate, 
$\beta_i$ denotes the transmission coefficient for strain~$i$,~$\gamma_i$
denotes the recovery rate from strain~$i$,~$\th_i$ is the rate at which
immunity to re-infection by strain~$i$ wanes and~$\th_{12}:=\th_3$ is the rate at
which immunity to re-infection by both strains wanes.  Finally,~$\si_{i} $
is the relative susceptibility to strain~$i$ for an individual previously
infected with and recovered from strain~$j$ ($i\ne j$), so that~$\si_{i}
=0$ corresponds to total cross-immunity,~$0<\si_{i} <1$ corresponds to
reduced susceptibility (partial cross-immunity) and~$\si_{i} >1$
corresponds to enhanced susceptibility, called ADE.
\beD[ADE] \cite{ferguson1999effect,schwartz2005chaotic, cummings2005}. For a multi-strain model, we say a secondary infection parameter
$\si_{1}$  has Antibody-Dependent Enhancement (ADE) if it is bigger than $1$.
\eeD

\beA In this paper, we consider only the mass action model obtained when $\al_i=0$.
\eeA
Then:
\BEN

\im The two forces of infection acting on $s$ are:
\begin{equation}
F_i=\beta_i {J_i}, \quad J_i=i_i+ \eta_i  y_i,  i=1,2,
\end{equation}

\im and the  forces of infection acting on
$r_i, i=1,2 $,
are: \begin{equation}
\si_{1} \beta_1 {J_1 }:=\T \beta_1 {J_1 }, \si_{2} \beta_2 {J_2  }:=\T \beta_2 {J_2  },
\end{equation}
where $\si_{1},\si_{2}$ denote   decrease  or increase factors of the susceptibility to secondary
infections.


  \EEN
{\bf Assumptions on model parameters}:
The feasible range of the problem parameters is
\begin{subequations}\label{eq:condParms}
\begin{equation}
\beta_i>0,\quad \gamma_i>0,\quad \mu\ge0,\quad \th_k\ge0,\quad
\eta_i> 0,\quad \sigma_{i}\ge0,\qquad i=1,2, \quad k=1,2,3.
\end{equation}
Assume that the susceptible group is replenished by demographic turnover ($\mu>0$) and/or by waning of the immune response generated following infections~($\th_k>0$),
\begin{equation}
\max\{\mu,\th_1,\th_2,\th_{12}\}>0
\end{equation}
\end{subequations}
(this enables the system to converge to an endemic equilibrium, rather than gradually exhausting the susceptible pool and converging to a disease-free state).

\beR
Some interesting features of the model, already revealed in \cite{nuno2005dynamics,Chung}, are: \BEN \im The expected form
$$R_1 >1,\Rot <1,$$ of the interior of the local stability domain of the boundary fixed point  $E_1$ (with an analog statement for $E_2$),  may be proved by factoring the full Jacobian
at $E_1$.
\im An  interior \fp\  is believed to exist, and may be proved under additional conditions -- see \cite[Thm 2.2]{GavCEP}. The general case  is a hard problem, even in the particular case of Chung \& Lui.

\im It appears from simulations that   when both \fbp s exist and are  unstable, either  a unique endemic point,  or an attracting periodic cycle may exist. The separation between the two is not fully understood (note that the Jacobian
at the  endemic point does not factor, and thus the respective Hopf bifurcation  occurs in dimension 7).

\im   Existence of the endemic point implies LAS --  see \cite[Thm 1.1, 1.2]{Chung}.

\EEN
\eeR
\beO[factoring the full Jacobian
at  boundary fixed points]
Is factoring the full Jacobian
at  boundary fixed points always possible for some class of multi-strain  models?
\eeO

\subsection{Getting the minimal siphons,  the DFE, the ngm, and  the reproduction functions with \texttt{bdAn}}
After inputting the model as a pair (RN,rts) we initiate the analysis by
calling the simplest bdAn (boundary analysis) module:
\begin{verbatim}
{RHS, var, par, cp, mSi, Jx, Jy, E0, K, R0A, ngm, infV} =
  bdAn[RN, rts];
Print["RHS=", RHS // FullSimplify // MatrixForm, "mSi=", mSi, " K= ",
K // MatrixForm, infV];
  \end{verbatim}

  The outputs of bdAn are:
  \begin{enumerate}
\item \textbf{RHS}: Right-hand side vector of the ODE system
\item \textbf{var}: List of all variables
\item \textbf{par}: List of all \para s
\item \textbf{cp}: List of positivity constraints for all \para s
\item \textbf{mSi}: Minimal siphons as variables (obtained by \texttt{minSiph})
    \item \textbf{Jx}: Jacobian of infection/invasion \var s (obtained by \texttt{NGM})
        \item \textbf{Jy}: Jacobian of non-infection \var s (obtained by \texttt{NGM})
            \item \textbf{E0}: DFE condition
                \item \textbf{K}: \NGM\ (obtained by \texttt{NGM})
                    \item \textbf{R0A}: List of non-zero \eval s of K
                     \item \textbf{ngm}: complete output of \texttt{NGM}, including $F,V$, and the alternative $K_d$ \NGM\
                         \item \textbf{infV}: the order of inf variables used in \texttt{NGM}.

\end{enumerate}

   This  reveals that the minimal siphons  correspond precisely to the two infectious strains: $i_1,i_{21}$  and $i_2,i_{12}$.
    The total  siphon is the union of the two minimal siphons, $\x=(i_1,i_{21},i_2,i_{12})$, $\y=(s,r_1,r_2,R=r_{12})$ .

At the DFE $E_0$, we find all species are $0,$ except $\sd=\fr{\La}{\mu}$.

The NGM $K$ is outputted in the order $i_1,i_2,i_{21},i_{12}$, and  after permutation to $i_1,i_{21},i_2,i_{12}$,  has a block structure $K=\left(
\begin{array}{cccc}
 \frac{\beta _1 s}{\gamma _1+\mu } & \frac{\beta _1 \eta _1 s}{\gamma _1+\mu } & 0 & 0 \\
 \frac{\beta _1 r_2 \sigma _1}{\gamma _1+\mu } & \frac{\beta _1 \eta _1 r_2 \sigma _1}{\gamma _1+\mu } & 0 & 0 \\
 0 & 0 & \frac{\beta _2 s}{\gamma _2+\mu } & \frac{\beta _2 \eta _2 s}{\gamma _2+\mu } \\
 0 & 0 & \frac{\beta _2 r_1 \sigma _2}{\gamma _2+\mu } & \frac{\beta _2 \eta _2 r_1 \sigma _2}{\gamma _2+\mu } \\
\end{array}
\right)$.
\beR Note the diagonal structure of the NGM, which could be takes as definition
of \mSt s, and explains the typical ``max" formulas that appear in multi-strain GLV or multi-strain SIR.

\eeR

The two diagonal blocks are the NGMs of the boundary \fp s defined
by the two minimal siphons. They have one $0$ \eig, and the reproduction functions (the unique positive \eig s, in this case), are precisely the traces of the two blocks, given by $$\bc R_1(\y)=\frac{\beta _1 \left(\eta _1 r_2 \sigma _1+s\right)}{\gamma _1+\mu },\\
R_2(\y)=\frac{\beta _2 \left(\eta _2 r_1 \sigma _2+s\right)}{\gamma _2+\mu }.\ec$$

It may be checked, both by direct stability analysis, and by the NGM method, that:
\beL The DFE, with $\sd=\fr{\La}{\mu}$, and all other \coo s zero, is unstable if
\be{R0} R_0=\max [R_1 \pr{s_0},R_2 \pr{s_0}]=\max \left[\fr{\sd}{s_1},\fr{\sd}{s_2}\right]>1.\ee
\eeL

\beR The \repF\ are increasing, and the DFE instability inequalities stated in Open Problem \ref{o:ins} hold.  Hence, the invasion stability conditions divide the parameter space in 4 pieces.
\eeR

\beR There is a version of bdAn, bdCo, which assumes that the all boundary points are rational, and furnishes them, and we could have used it here.
\How we have preferred to give the more gradual analysis, where the non-DFE points are analyzed  in a second stage.
\eeR
\subsection{Obtaining the 
single-strain endemic equilibrium points with \texttt{bdFp} and \texttt{inv2}}\label{s:bendG}
Put now  \be{eq:Ri}
\mathcal{R}_i:=\frac{\beta _i}{\gamma _i+\mu }, R_i:=R_i(E_0)= \mathcal{R}_i \sd.\ee

Assume \fno\ \wlo\ that $E_1$ is the more virulent strain,  (excluding the non-generic equality case), \ie\ $\mR_{2} < \mR_{1}  \Eq R_2 <  R_1$.

The  single-strain endemic equilibrium points of~\eqref{eq:model} are found by solving the \fp\ system under the siphon conditions, using \texttt{bdFp}.
For each boundary system, the output consists of a pair consisting of all rational solutions, and of a polynomial \satd\ by the non-rational solutions.
Here, the output of \texttt{bdFp}
\begin{verbatim}
bdfp = bdFp[RHS, var, mSi];
Print["rat sols on first siphon facet are"]
bd1 = bdfp[[1, 1]] // FullSimplify
bdfp[[1, 2]]
\end{verbatim}
 reveals that all solutions are rational, and that on each of the two siphon facets we have three \fp s, one of which is the DFE, and one of which is not \nne. The indexes of the correct solution (2,2)
 are then passed on to
 \texttt{invN2}, which computes the invasion numbers and a numerical example
 under which the system is conjectured to be persistent:
\begin{verbatim}
{E1, E2, R12, R21, coP} =
  invN2[bdfp[[1, 1]], bdfp[[2, 1]], R0A, E0, par, cp, 2, 2];
Print["invasion numbers R12, R21 are ", R12 // Apart, R21 // Apart]
\end{verbatim}

We may conclude that:
\beT[stability of the boundary fixed points for the Gavish model]\label{thm:EE1}
Assume the parameters~$\mu$, $\{\gamma_i\}_{i=1}^2$, $\{\th_k\}_{k=1}^3$, $\{\eta_i\}_{i=1}^2$,~$\si_{2}$ and~$\sigma_{1}$  satisfy~\eqref{eq:condParms}. Then:

 1) If $R_1>  1 \Eq \sd>s_1$, then,
the system~\eqref{eq:model} has a unique single-strain endemic equilibrium~$\phi^{E,1}$ with~$i_2=i_{12}=0$ and~$i_1>0$. Putting $$\T{\ga}_1=\frac{\ga_1}{\ga_1+\th_1+\mu },$$ this solution satisfies
\bea
&&\Scale[.9]{s_1=\frac1{\mR_1},\quad r_1^{E,1}=\T{\ga}_1\pr{\sd-s_1}=\T{\ga}_1 \sd\pr{1-\fr{1}{R_1}},
i_1^{E,1}=(1-\T{\ga}_1)\pr{\sd-s_1}=(1-\T{\ga}_1)\sd\pr{1-\fr{1}{R_1}}}\\&&
\quad i_2^{E,1}=i_{12}^{E,1}=r_2^{E,1}=i_{21}^{E,1}=0,\\
&&s_1+i_1^{E,1}+r_1^{E,1}=\sd,
\eea
and an analog statement holds for a unique single-strain endemic equilibrium~$\phi^{E,2}$ with~$i_1=i_{21}=0$ and~$i_2>0$.

2) The invasion number $\Rot$ is given by:
\be{R12G}\Rot=\mR _2 \left(\sigma _2 \eta _2 r_1^{E,1} +s_1\right)
=\mR _2 \left(\sigma _2 \eta _2 \T{\ga}_1 \sd\pr{1-\fr{1}{\mR_1}} +\frac1{\mR_1}\right)\ee
see \cite[5]{GavCEP}.

3) $\phi^{E,1}$ is LAS if \( \Rot <1\)  and unstable if \( \Rot >1\).

\eeT

\Prf  For the last point, we exploit the key observation
made by \cite{nuno2005dynamics,Chung} (for the simpler particular case of permanent immunity) that the full Jacobian
at $E_1$ factors into one quadratic with negative roots when $r_1>1$, and linear terms, only one of which may have positive roots. The condition that this root \cite[(2.2)]{Chung} is positive may be finally be expressed as
\( \Rot >1\).

\QED

\beR The  proof offered  in \cite[Thm (2.1)]{GavCEP}    is incomplete, since the conditions of \cite{Van} are not checked, and this is not at all trivial at the non-DFE \bfp s; however, as implemented in our .nb file, the  characteristic polynomial factorization proof of \cite{Chung} is easily checked to extend to this case. \eeR

\subsection{The endemic point, in the particular case of \cite{Chung}}

The existence of the  endemic point is a hard problem, still unresolved even in the particular case of Chung \& Lui \cite{Chung}, when $\eta_1=\eta_2=1, \th_i=0, i=1,2,3, \La=\mu$.

For this model studied already in \cite{Chung}, the \fp\ system may be reduced to the scalar equation
$a_0 + a_1 s+a_2 s^2+a_{12} s^3,$ where
$$\bc a_0=\sigma _1 \sigma _2 (\gamma _1+\mu ){}^2 (\gamma _2+\mu ){}^2\\
a_1=(\gamma _1+\mu ) (\gamma _2+\mu ) (\sigma
_1 (\gamma _1+\mu ) (\gamma _2+\mu )-\sigma
_2 (\sigma _1 (\gamma _2+\mu ) (\beta _1+\gamma
_1+\mu )  \\  ~~~~~~~~~+\beta _2 \sigma _1 (\gamma _1+\mu )-((\gamma _1+\mu ) (\gamma _2+\mu ))))
\\a_2=(-\beta _2) \sigma _1 (\gamma _2+\mu ) (\gamma _1+\mu ){}^2+\beta _1
(\gamma _2+\mu ) (\gamma _1+\mu ) ((\sigma _1-1) \sigma _2 (\gamma _2+\mu
)-\mu  \sigma _1) \\ ~~~~~~~~~+\beta _2 \sigma _2 (\sigma _1 (\gamma _1 \mu  (\beta _1+2 (\gamma _2+\mu ))+\mu  (\beta _1+\mu ) (\gamma _2+\mu )+\gamma _1^2 (\gamma _2+\mu ))-\mu  (\gamma _1+\mu ) (\gamma _2+\mu ))
       \\a_{12}=\beta _1 \beta _2 \mu  (\sigma _1 (\gamma _1+\mu )+\sigma _2 (-\sigma _1 (\gamma _1+\gamma _2+\mu )+\gamma _2+\mu ))\ec,$$
       which confirms  \cite[(6b)]{Gavnew}.
A generalization including the immunity \para s may be found in the .nb file.

Here is the plot of an instance of the partition. Further experiments are necessary to  determine if Hopf bifurcations are possible (in this instance there is none).

\figu{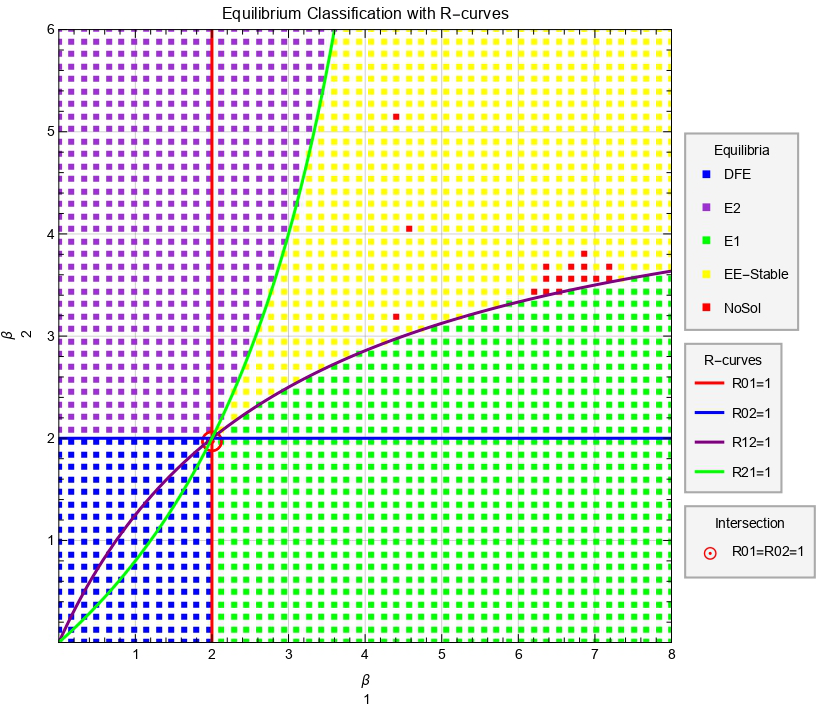}{The stability partition scan obtained numerically in GavScan.nb almost coincides with the symbolic one obtained using the explicit reproduction and invasion numbers. The file GavScan.nb may be modified (simply change the reactions) to obtain partition scans in cases where the symbolic picture is incomplete (for example when Hopf bifurcations are possible)}{.7}

\beO[coexistence of several locally stable \fp s] \lbl{coe} For general two-strain models,
is it possible that  both strains coexist  and are locally stable? (recall this is impossible for LVK models)
\eeO

We end this section with a fundamental open problem.
\beO [The persistence conjecture] \lbl{o:per}The  examples studied here and others  suggest that for a positive ODE, the instability of all fixed boundary points (or maybe even only the instability of the disease-free equilibrium and  the boundary equilibria whose zero set is a maximal non-DFE siphon) ensures the persistence.
\eeO

For some positive results in particular cases, see \cite{AdLS11,PantGAS,CraNazPan,GopMilShiu,
craciun2019polynomial}.

For another particular case where the persistence conjecture holds, we recall  the following:
\beP (\cite[Thm 4.19]{Bulh})
Two-strain models which live on a compact \nne\ invariant subset and have   $R_j > 1, j =1, 2$,
 and $R_{\T 2}^1> 1, R_{\T 1}^2>1,$ so that both strains exist, but are unstable, are uniformly persistent.
 \eeP

\ssec{Can the persistence and global stability of this two strain model 
be resolved by removing intermediates?
}

This question is inspired by the CRN papers \cite{FFW,FWF}, which provide conditions under which the removal of  intermediates does not affect persistence. This suggests that there might be a relation between the  persistence or GAS property of this model when $r_0>1$, and those of   the model where the intermediate states $r_i$ are removed.  

 \How the \cite{FFW,FWF} theory does not apply:  the  cross-immunity assumption $\si_i>0$ implies that $r_1, r_2$ are essential species, which do not \saty\ the FFW reducibility assumption that intermediates may not appear in mixed complexes
which mix ``intermediates" and ``non-intermediates" (r1+i2, r1+i12, r2+i1, r2+i21).
 Cross-immunity creates a ``biological memory" that prevents mathematical reduction; the two-strain  model with cross-immunity is ``FFW-irreducible".

\section{A   SI$^2$R  model with co-infections of  \cite{GK}, which has more reproduction functions than minimal siphons;SI2Coinf.nb}\label{s:Koz}

Co-infection models \cite{Regoes2004CTL,Vautrin2007WolbachiaMultipleInfections,Althaus2008Epitope,
Lipsitch2009NoCoexistence,
Ganusov2011Escape,
Alizon2013Coinfection,Lion2013MultipleInfections,
Gao2016CoinfectionTwoDiseases,Mitchell,Belluccini} pose a  mathematical dilemma in our construction. Consider  for  example  the logistic  growth SI$^2$R model  of  \cite{GK}, shown in  figure
\ref{fig:coinfection_flow} below;  should it be viewed as two-strain, or three strain?

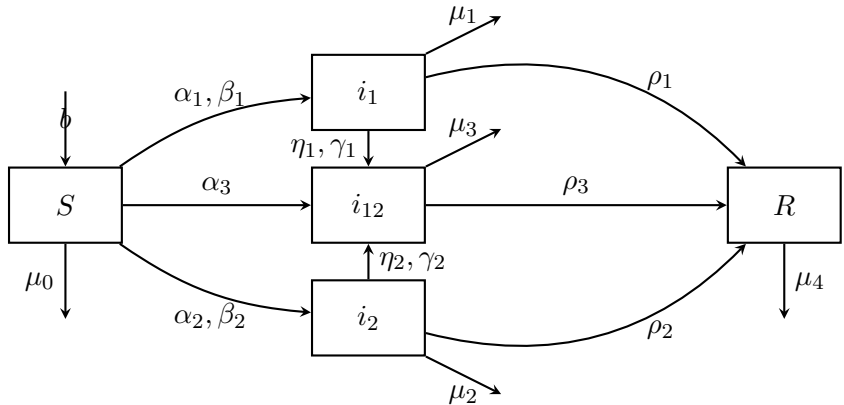
\begin{figure}[H]
\centering
\begin{tikzpicture}[
  node distance=2.5cm,
  compartment/.style={
    draw,
    rectangle,
    minimum height=1cm,
    minimum width=1.5cm,
    thick,
    align=center
  },
  arrow/.style={
    ->,
    thick,
    >=stealth
  }
]

\node[compartment] (S) {$S$};
\node[compartment] (I1) [right=of S, yshift=1.5cm] {$i_1$};
\node[compartment] (I12) [right=of S] {$i_{12}$};
\node[compartment] (I2) [right=of S, yshift=-1.5cm] {$i_2$};
\node[compartment] (R) [right=of I12, xshift=1.5cm] {$R$};

\draw[arrow] (S) to[bend left=15] node[above, midway, draw=none, fill=none] {$\alpha_1, \beta_1$} (I1);
\draw[arrow] (S) to[bend right=15] node[below, midway, draw=none, fill=none] {$\alpha_2, \beta_2$} (I2);

\draw[arrow] (S) to node[above, draw=none, fill=none] {$\alpha_3$} (I12);

\draw[arrow] (I1) to node[left, draw=none, fill=none] {$\eta_1, \gamma_1$} (I12);
\draw[arrow] (I2) to node[right, draw=none, fill=none] {$\eta_2, \gamma_2$} (I12);

\draw[arrow] (I1) to[bend left=30] node[above, pos=0.7, draw=none, fill=none] {$\rho_1$} (R);
\draw[arrow] (I12) to node[above, draw=none, fill=none] {$\rho_3$} (R);
\draw[arrow] (I2) to[bend right=30] node[below, pos=0.7, draw=none, fill=none] {$\rho_2$} (R);

\draw[arrow] ([yshift=1cm]S.north) -- (S.north)
    node[midway,yshift=4pt] {$b$};

\draw[arrow] (S.south) -- ++(0,-1) node[midway,left] {$\mu_0$};
\draw[arrow] (I1.north east) -- ++(1,0.5) node[midway,above] {$\mu_1$};
\draw[arrow] (I2.south east) -- ++(1,-0.5) node[midway,below] {$\mu_2$};
\draw[arrow] (I12.north east) -- ++(1,0.5) node[midway,above] {$\mu_3$};
\draw[arrow] (R.south) -- ++(0,-1) node[midway,right] {$\mu_4$};

\end{tikzpicture}
\caption{Flow diagram for \cite{GK} two-strain coinfection model}
\label{fig:coinfection_flow}
\end{figure}

where the following notation was used:
\begin{itemize}
\item $i_1$ and $i_2$ are infected classes from strain 1 and strain 2 respectively.
\item $i_{12}=i_3$ is the compartment consisting of individuals infected by both diseases.
\item $R$ represents the recovered class.
\end{itemize}

Due to the permanent immunity assumed,  we may remove R  from the analysis;  the total removal rates of the infected will be denoted by $\mu_i=\mu_i'+ \rho_i, i=1,2,3$.

We open now a parenthesis about the similar Volpert/Petri/\DSR\ (DSR) bipartite graph, which plays a big role in CRNT.  The graph below,  for  pairs of species connected by more than one reaction in Figure \ref{fig:coinfection_flow}, includes an edge for  each  reaction, and it indicates the complete rate  inside  a circle, which represents the reaction. However,   for  pairs of species connected by one reaction only, it keeps the structure of \ref{fig:coinfection_flow}, except that is indicates the complete rate, not just its constant. The DSR graph is obtained by trasforming all edges in an edge with a circle in the middle. Clearly, the DSR representation is more complete than the flow diagram in Figure \ref{fig:coinfection_flow}.

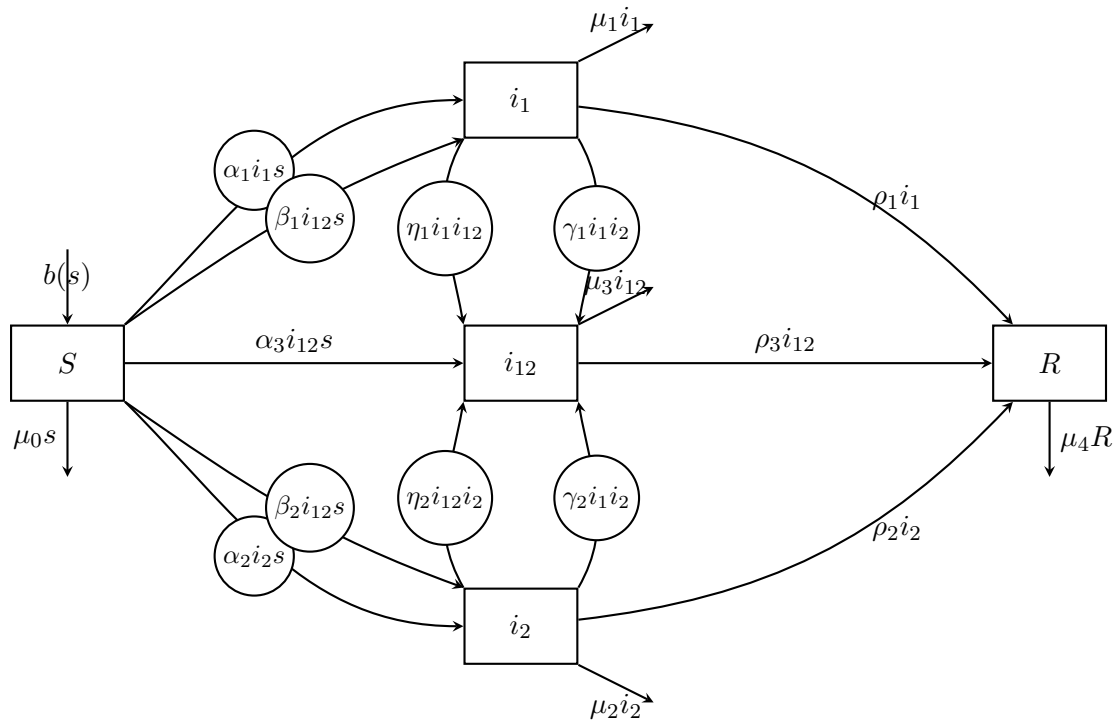
\begin{figure}[H]
\centering
\begin{tikzpicture}[
  node distance=2.5cm,
  compartment/.style={
    draw,
    rectangle,
    minimum height=1cm,
    minimum width=1.5cm,
    thick,
    align=center
  },
  arrow/.style={
    ->,
    thick,
    >=stealth
  },
  rate/.style={
    draw,
    circle,
    fill=white,
    minimum size=0.8cm,
    inner sep=2pt,
    font=\small
  }
]
\node[compartment] (S) {$S$};
\node[compartment] (I1) [right=of S, xshift=2cm, yshift=3.5cm] {$i_1$};
\node[compartment] (I12) [right=of S, xshift=2cm] {$i_{12}$};
\node[compartment] (I2) [right=of S, xshift=2cm, yshift=-3.5cm] {$i_2$};
\node[compartment] (R) [right=of I12, xshift=3cm] {$R$};

\draw[arrow] (S.north east) to[out=45, in=180]
    node[pos=0.5, rate, xshift=-0.3cm] {$\alpha_1 i_1 s$} (I1.west);
\draw[arrow] (S.north east) to[out=35, in=200]
    node[pos=0.5, rate, xshift=0.3cm] {$\beta_1 i_{12} s$} (I1.south west);

\draw[arrow] (S) to node[above, draw=none, fill=none] {$\alpha_3 i_{12} s$} (I12);

\draw[arrow] (S.south east) to[out=-45, in=180]
    node[pos=0.5, rate, xshift=-0.3cm] {$\alpha_2 i_2 s$} (I2.west);
\draw[arrow] (S.south east) to[out=-35, in=160]
    node[pos=0.5, rate, xshift=0.3cm] {$\beta_2 i_{12} s$} (I2.north west);

\draw[arrow] (I1.south west) to[out=-120, in=100]
    node[pos=0.5, rate] {$\eta_1 i_1 i_{12}$} (I12.north west);
\draw[arrow] (I1.south east) to[out=-60, in=80]
    node[pos=0.5, rate] {$\gamma_1 i_1 i_2$} (I12.north east);

\draw[arrow] (I2.north west) to[out=120, in=-100]
    node[pos=0.5, rate] {$\eta_2 i_{12} i_2$} (I12.south west);
\draw[arrow] (I2.north east) to[out=60, in=-80]
    node[pos=0.5, rate] {$\gamma_2 i_1 i_2$} (I12.south east);

\draw[arrow] (I1) to[bend left=20] node[above, pos=0.7, draw=none, fill=none] {$\rho_1 i_1$} (R);
\draw[arrow] (I12) to node[above, draw=none, fill=none] {$\rho_3 i_{12}$} (R);
\draw[arrow] (I2) to[bend right=20] node[below, pos=0.7, draw=none, fill=none] {$\rho_2 i_2$} (R);

\draw[arrow] ([yshift=1cm]S.north) -- (S.north)
    node[midway,yshift=4pt] {$b(s)$};

\draw[arrow] (S.south) -- ++(0,-1) node[midway,left] {$\mu_0 s$};
\draw[arrow] (I1.north east) -- ++(1,0.5) node[midway,above] {$\mu_1 i_1$};
\draw[arrow] (I2.south east) -- ++(1,-0.5) node[midway,below] {$\mu_2 i_2$};
\draw[arrow] (I12.north east) -- ++(1,0.5) node[midway,above] {$\mu_3 i_{12}$};
\draw[arrow] (R.south) -- ++(0,-1) node[midway,right] {$\mu_4 R$};
\end{tikzpicture}
\caption{Flow diagram for \cite{GK} two-strain coinfection model in DSR format}
\label{fig:DSRKoz}
\end{figure}

For a historical note on these  Volpert/Petri/DSR graphs, see Table  \ref{t:histNoteVolpertgraphs}.

\begin{table}[h!]
\centering
\caption{Comparison of common graph representations of chemical reaction networks.}
\label{t:histNoteVolpertgraphs}
\resizebox{\textwidth}{!}{%
\begin{tabular}{@{}lllll@{}}
\toprule
\textbf{Name} & \textbf{Origin} & \textbf{Nodes} & \textbf{Edge Direction} & \textbf{Typical Use} \\ \midrule
Volpert graph & Chemical kinetics (Volpert, 1960s) & Species and reactions & Directed & Reachability, persistence \\
Petri net & Computer science (Petri, 1962) & Places and transitions & Directed & Discrete simulation, invariants \\
SR graph & CRN theory (Craciun, Feinberg, 2000s) & Species and reactions & Undirected (signed) & Injectivity tests \\
DSR graph & Dynamical systems (Banaji \& Craciun, 2010s) & Species and reactions & Directed (signed) & Stability, monotonicity \\ \bottomrule
\end{tabular}%
}
\end{table}

The reaction network representation of the  \cite{GK} model is:
$$\left(
\begin{array}{cc}
 0 \to \text{s} & b(s) \\
 \text{i1}+\text{s}\to 2 \text{i1} & \alpha _1 i_1 s \\
 \text{i2}+\text{s}\to 2 \text{i2} & \alpha _2 i_2 s \\
 \text{i12}+\text{s}\to 2 \text{i12} & \alpha _3 i_3 s \\
 \text{i1}+\text{i12}\to 2 \text{i12} & \eta _1 i_1 i_3 \\
 \text{i12}+\text{i2}\to 2 \text{i12} & \eta _2 i_3 i_2 \\
 \text{i1}+\text{i2}\to \text{i12}+\text{i2} & \gamma _1 i_1 i_2 \\
 \text{i1}+\text{i2}\to \text{i12}+\text{i1} & \gamma _2 i_1 i_2 \\
 \text{i12}+\text{s}\to \text{i12}+\text{i1} & \beta _1 i_3 s \\
 \text{i12}+\text{s}\to \text{i12}+\text{i2} & \beta _2 i_3 s \\
 \text{s}\to 0 & \mu _0 s \\
 \text{i1}\to 0 & i_1 \mu _1 \\
 \text{i2}\to 0 & i_2 \mu _2 \\
 \text{i12}\to 0 & i_3 \mu _3 \\
\end{array}
\right).
$$

The EpidCRN  script {\emph bdAn}, furnishes many preliminary results for an epidemiologic analysis:
\begin{verbatim}
{RHS, var, par, cp, mSi, Jx, Jy, E0, K, R0A, infVars, gam, ng} =
  bdAn[RN, rts];
  {edg, cyc, graph} = IGMS[RN, mSi];
Print["siphons are", mSi, " edges are",
 edg, "  DFE is", E0, " repr. functions=  ", R0A]; F =
 ng[[2]]; V = ng[[3]];
Print["K=", K // MatrixForm, "F=", F // MatrixForm,
"V=",V // MatrixForm];
\end{verbatim}
  {\emph bdAn} provides:
  \BEN \im The RHS of the ODE.
 \im  The minimal siphons mSi, which, assuming $\beta>0, \gamma>0$, are \{i1,i12\},\{i2,i12\};    the total DFE siphon is therefore \{i1,i2,i12\}.

\im The subsequent command {\emph IGMS[RN,mSi]} reveals that the IGMS is a cycle
\figu{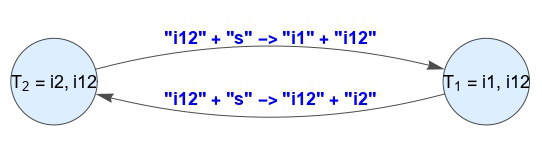}{IGMS for  \cite{GK} two-strain  model; edges denote
  coinfection reactions}{.9}
\im
The stability of the jacobian $J_y=- \left(i_3 \left(\alpha _3+\beta _1+\beta _2\right)+\alpha _1 i_1+\alpha _2 i_2+\mu _0\right)<0$ is one of the necessary
conditions we ask from a ME model.
\im Another one is the existence of a \regS\ of $J_x$, or of a matrix $M$ appearing in the factorization of the  \fp\ equations \mbw:
\be{submodel2}\bff 0=\bep
b-s \left(i_3 \left(\alpha _3+\beta _1+\beta _2\right)+\alpha _1 i_1+\alpha _2 i_2+\mu _0\right)\\
\bep \alpha _1 s-\gamma _1 i_2-\eta _1 i_3-\mu _1&0& \beta _1  s \\
 0&\alpha _2 s-\gamma _2 i_1-\eta _2 i_3-\mu _2&\beta _2  s \\
 \eta _1 i_3+\gamma _1 i_2&\eta _2 i_3 +\gamma _2 i_1 &\alpha _3 s-\mu _3\eep

\bep i_1\\i_2\\i_3
\eep
\eep
\ee
\bea =\bep
b-s \left(i_3 \left(\alpha _3+\beta _1+\beta _2\right)+\alpha _1 i_1+\alpha _2 i_2+\mu _0\right)\\
\bep \alpha _1 s-\gamma _1 i_2-\eta _1 i_3-\mu _1&0& \beta _1  s \\
 0&\alpha _2 s-\gamma _2 i_1-\eta _2 i_3-\mu _2&\beta _2  s \\
 \gamma _1 i_2&\gamma _2 i_1 &\eta _1 i_1+\eta _2 i_2 +\alpha _3 s-\mu _3\eep

\bep i_1\\i_2\\i_3
\eep
\eep
\eea

Denoting by $M$ the matrix multiplied by $\bep i_1\\i_2\\i_3
\eep$ in the first \rep\ in \eqr{submodel2} above, and spitting $M=F-V$ as the difference of a \nne\ matrix and of a matrix with positive inverse
$$M=s\bep \alpha_1 & 0& \beta_1 \\
0&\alpha_2 &\beta_2 \\
0&
0& \alpha_3   \eep+
\bep -\gamma_1 i_2 -\eta_1 i_3  -\mu_1& 0& 0\\
0&-\gamma _2 i_1-\eta_2 i_3 -\mu_2&0\\
\gamma_1 i_2+ \eta_1 i_3&
\gamma _2 i_1+\eta_2 i_3&   -\mu_3 \eep
$$
we may conclude that our model \sats\ our definition of ME models.
\im The second \rep\ in \eqr{submodel2} reveals that when $\beta=\ga=0$, this is a Lotka-Volterra-Kolmogorov model.
\im
The next generation matrix may be made triangular \  $$K=F V^{-1}=\left(
\begin{array}{ccc}
 \frac{\alpha _1 s}{\mu _1} & 0 & \frac{\beta _1 s}{\mu _3} \\
 0 & \frac{\alpha _2 s}{\mu _2} & \frac{\beta _2 s}{\mu _3} \\
 0 & 0 & \frac{\alpha _3 s}{\mu _3} \\
\end{array}
\right)$$
by ordering the variables as $i_1,i_2,i_{12}$. This suggests that the most convenient ordering of the DFE support must start with elements which appear only in one siphon.

Finally, triangularity reveals that $R_0=s_0 \max_{i=1}^3 [\mR_i],$ where
\be{mRK} \mR_i=\frac{\alpha_i}{\mu_i},\quad i=1,2,3,\ee
and where  $s_0=\fr{b}{\mu_0}$ in the linear growth case (and equal to the carrying capacity in the logistic case).
This concludes the stability analysis of the DFE.
\im
 ngm is the list of all  the outputs of NGM, which includes the matrices F,V.

\im Using the notation \eqr{mRK},
we may express some of the boundary fixed points which have precisely one positive infection  variable in a revealing form:
$$E_1=(s= \frac{1}{\mR _1},i_1= \fr{\mu _0}{\alpha _1 } \pr{\mR _1 s_0- 1},i_2= 0,i_3 =0),$$
with   similar formulas for $E_2, E_3$. These points exist iff
\be{RK}R_i:=\mR _i s_0 >1, i=1,2,3.\ee

 \be{RfK}R_i(s):= \mR _i s,\ee
  called  reproduction functions and applied at the s-value of the other fixed points (besides \eqr{RK}), will be seen  to intervene in the stability analysis in Section \ref{s:Koa}.

\EEN

\subsection{Further analytical results}\lbl{s:Koa}
Recall the  reaction network from Section \ref{s:Koz}:
$$\left(
\begin{array}{cc}
 0 \to \text{s} & b(s) \\
 \text{i1}+\text{s}\to 2 \text{i1} & \alpha _1 i_1 s \\
 \text{i2}+\text{s}\to 2 \text{i2} & \alpha _2 i_2 s \\
 \text{i12}+\text{s}\to 2 \text{i12} & \alpha _3 i_3 s \\
 \text{i1}+\text{i12}\to 2 \text{i12} & \eta _1 i_1 i_3 \\
 \text{i12}+\text{i2}\to 2 \text{i12} & \eta _2 i_3 i_2 \\
 \text{i1}+\text{i2}\to \text{i12}+\text{i2} & \gamma _1 i_1 i_2 \\
 \text{i1}+\text{i2}\to \text{i1}+\text{i12} & \gamma _2 i_1 i_2 \\
 \text{i12}+\text{s}\to \text{i1}+\text{i12} & \beta _1 i_3 s \\
 \text{i12}+\text{s}\to \text{i12}+\text{i2} & \beta _2 i_3 s \\
 \text{s}\to 0 & \mu _0 s \\
 \text{i1}\to 0 & i_1 \mu _1 \\
 \text{i2}\to 0 & i_2 \mu _2 \\
 \text{i12}\to 0 & i_3 \mu _3 \\
\end{array}
\right),
$$
with ODE
\bea X'=\left(
\begin{array}{c}
 b-s \left(i_3 \left(\alpha _3+\beta _1+\beta _2\right)+\alpha _1 i_1+\alpha _2 i_2+\mu _0\right) \\
 i_1 \left(-\gamma _1 i_2-\eta _1 i_3-\mu _1+\alpha _1 s\right)+\beta _1 i_3 s \\
 i_2 \left(-\gamma _2 i_1-\eta _2 i_3-\mu _2+\alpha _2 s\right)+\beta _2 i_3 s \\
 \left(\gamma _1+\gamma _2\right) i_1 i_2+i_3 \left(\eta _1 i_1+\eta _2 i_2-\mu _3+\alpha _3 s\right) \\
\end{array}
\right)
\eea

The   parameters of the system are:
\begin{itemize}
\item $b$ is the intrinsic growth of $s$. \cite{GK} study the logistic case with $b=r*s*(1 - s/K)$, containing two further \para s, the per capita birth rate $r$,
and the carrying capacity
$K$, but we study  the case when $b$ is  constant.
\item
$\gamma_i$ is the rate at which individuals infected with one strain get infected with  the other strain after a meeting with single infected individuals, and move to the coinfected class ($i=1,2$).
\item $\mu_i',i=1,2,3 $ (not shown in figure) are the death rates of the infected \com  s.
\item
$\alpha_1$, $\alpha_2$, $\alpha_3$  are the rates of transmission of strain 1, strain 2 and both strains (\ie\ coinfection),
\item $\beta_i$ is the rate at which susceptibles contract disease $i=1,2$ but not the other disease, from a coinfected individuals.

\item
$\eta_i, i=1,2$ is the rate at which individuals infected by one strain gets coinfected, after a meeting  with a coinfected individual.
\end{itemize}

The following vector notations will  be used:
$$
\alpha=(\alpha_1,\alpha_2,\alpha_3),\;\;\mu=(\mu_0,\mu_1,\mu_2,\mu_3),\;\;\eta=(\eta_1,\eta_2),
\;\;\mbox{and}\;\;
$$

The parameters are split into three groups:
\begin{equation}\label{Bparam}
{\mathcal B}=(b,\alpha,\mu,\eta), \text{ (basic) },\;\beta=(\beta_1,\beta_2), \;\;\;\gamma=(\gamma_1,\gamma_2), \; \; {\Omega }=(\gamma_1,\gamma_2,\beta_1,\beta_2)\; \text{ (coinfection).}
\end{equation}

     It is reasonable to assume that the coinfection constants from    ${\Omega }$, involved in terms describing double infection,   are small with respect to basic constants.

  Assume  that ${\mathcal B}$ are positive, and $\beta, \gamma $ are \nne.

   \beR
One interesting aspect of this model is the sensitivity of the IGMS structure to the \para s.  When $\beta>0, \gamma>0$, the minimal siphons are \{i1,i12\},\{i2,i12\}. When $\beta=0, \gamma>0,$ the unique minimal siphons is \{i12\}.

And when $\beta=0= \gamma$, there are three minimal siphons    \{\{i1\},\{i2\},\{i12\}\},  and the IGMS is a union of two cycles:
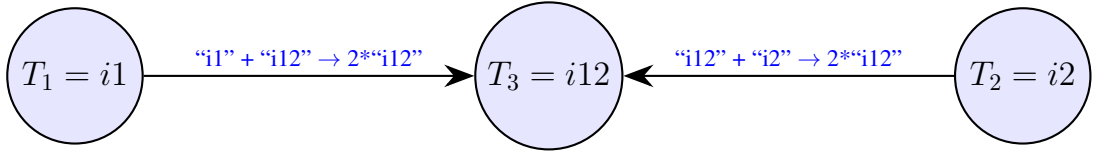
\begin{figure}[h!]
\centering
\begin{tikzpicture}[
  scale=0.9,  
  transform shape,  
  strain/.style={draw, circle, minimum size=2cm, thick, fill=blue!10, font=\Large},
  arrow/.style={-{Stealth[length=4mm]}, thick, >=stealth},
  reaction/.style={font=\small, color=blue}
]

\node[strain] (T1) at (0,0) {$T_1 = i1$};
\node[strain] (T3) at (7,0) {$T_3 = i12$};
\node[strain] (T2) at (14,0) {$T_2 = i2$};

\draw[arrow] (T1) -- node[above, reaction, pos=0.5] {``i1'' + ``i12'' $\to$ 2*``i12''} (T3);

\draw[arrow] (T2) -- node[above, reaction, pos=0.5] {``i12'' + ``i2'' $\to$ 2*``i12''} (T3);

\end{tikzpicture}
\caption{The IGMS of the $\beta=0= \gamma$ coinfection
Kozlov model   with coinfection strain i12 has two cycles. Strains i1 and i2 both lose to the coinfection strain i12 when they come into contact.}
\label{fig:coinfection_tree}
\end{figure}
\eeR

\Fr as noted already, when $\beta _1=\beta _2=\gamma _1=\gamma _2=0$, the infection part is a Lotka-Volterra-Kolmogorov model with antisymmetric interaction matrix, to be studied in the next section.

\ssec{The antisymmetric Lotka-Volterra-Kolmogorov case of the \cite{GK} model, with 0 coinfection \para s}\lbl{s:KLV}

With 0 coinfection \para s, $\beta_j=\ga_j=0, j=1,2,$ the \cite{GK} ODE has its infection part of  Lotka-Volterra-Kolmogorov type:
\bea X'=\left(
\begin{array}{c}
 \mu _0 s_0-s \left(i_3 \mu_3 \mR_3+\mu_1 \mR_1 i_1+\mu_2 \mR_2 i_2+\mu _0\right) \\
 i_1 \left(-\eta _1 i_3-\mu _1+\mu_1 \mR_1 s\right)\\
 i_2 \left(-\eta _2 i_3-\mu _2+\mu_2 \mR_2 s\right)\\
i_3 \left(\eta _1 i_1+\eta _2 i_2-\mu _3+\mu_3 \mR_3 s\right) \\
\end{array}
\right)
\eea

The infection equations \mbw:
\bea x'= \bep
 i_1 \\
 i_2 \\
i_3  \\
\eep \circ \pp{\bep \mu _1(s \mR_1-1)\\ \mu _2(s \mR _2 -1)\\ \mu _3(s \mR _3 -1)\eep + \bep0&0 &-\eta _1 \\
 0 &0&-\eta _2 \\
\eta _1 &\eta _2 &0
\eep .\bep
 i_1 \\
 i_2 \\
i_3  \\
\eep}:=x \circ (c+ A x)
\eea
with $A$ antisymmetric. The persistence - extinction dichotomy can thus be tackled via the theory of  Lotka-Volterra-Kolmogorov ODEs with anti-symmetric $A$ \cite{xu2024persistence}.

This model has six fixed points, all rational, given respectively by:
\begin{equation*}
\bc
E_1=\left\{s= \frac{1}{ \mR_1}:=s_1,
i_1= \frac{ \mu_0}{\mu _1}(s_0-s_1)= \frac{ \mu_0 s_1}{\mu _1}(s_0 \mR_1 -1):=\T i_1,i_2= 0,i_3= 0\right\},
\\E_2=\left\{s= \frac{1}{ \mR_2}:=s_2,i_1= 0,i_2=
\frac{\mu_0}{\mu _2}(s_0-s_2)= \frac{ \mu_0 s_2}{\mu _1}(s_0 \mR_2 -1):=\T i_2,i_3= 0\right\},\\
\Scale[.82]{E3=\left\{s= \frac{1}{ \mR_3}:=s_3,i_1= 0,i_2= 0,i_3=
\frac{\mu_0}{\mu _3}(s_0 -s_3):=\T i_3\right\},DFE=\left\{s= \frac{s_0 \mu_0}{\mu _0},i_1= 0,i_2= 0,i_3= 0\right\},}\\
\Scale[.8]{E13=\left\{s_{13}= \frac{s_0 \mu_0 \eta _1}
{\mu_1  \mu _3(\mR_1-\mR_3)+\eta _1 \mu _0},
i_1= \frac{\mu _3}{\eta _1}(1- \mR_3 s_{13})=\frac{\al _3}{\eta _1}(s_3 - s_{13}),i_3= \frac{\mu _1}{\eta _1}(\mR_1 s_{13}-1)= \frac{\al _1}{\eta _1}( s_{13}-s_1),i_2= 0\right\}},
\\E23=\left\{s_{23}= \frac{s_0 \mu_0 \eta _2}{\mu_2  \mu _3(\mR_2-\mR_3)+\eta _2 \mu _0},i_1= 0,
i_2=\frac{\mu _3}{\eta _2}(1-  \mR_3 s_{23}),i_3= \frac{\mu _2}{\eta _2}( \mR_2 s_{23}-1)\right\}
\ec
\end{equation*}

This reveals that the minimal siphon {i12=0} does not have any fixed point in its interior (however, it contains the two fixed points $E_1,E_2$ on the boundary).

Finally, the endemic point $EE$ has \coo s:
\bea &&
s= \frac{\eta _1 \mu _2-\eta _2 \mu _1}{\mu_2 \mR_2 \eta _1-\mu_1 \mR_1 \eta _2}:=\sE, i_3= \frac{\mu_1  \mu _2(\mR_1- \mR_2)}{\mu_2 \mR_2 \eta _1-\mu_1 \mR_1 \eta _2}= \frac{\mu_1  \mu _2(\mR_1- \mR_2) \sE}{\eta _1 \mu _2-\eta _2 \mu _1}\\&&
\Scale[1]{i_1= \frac{\mu _2 \left(\mu_3 \mR_3 \eta _2 \mu _1+\mu_2 \mR_2 \eta _1 \mu _3+\eta _1 \eta _2 \mu _0\right)+\mu_3 \mR_3 \eta _1 \left(-\mu _2^2\right)-\eta _2 \mu _1 \left(\mu_2 \mR_2 \mu _3+\eta _2 \mu _0\right)+s_0 \mu_0 \eta _2 \left(\mu_1 \mR_1 \eta _2-\mu_2 \mR_2 \eta _1\right)}{\left(\mu_2 \mR_2 \eta _1-\mu_1 \mR_1 \eta _2\right) \left(\eta _1 \mu _2-\eta _2 \mu _1\right)},}
\\&& \Scale[1]{i_2= \frac{\mu _1 \left(\mu_3 \mR_3 \eta _1 \mu _2+\mu_1 \mR_1 \eta _2 \mu _3+\eta _1 \eta _2 \mu _0\right)+\mu_3 \mR_3 \eta _2 \left(-\mu _1^2\right)-\eta _1 \mu _2 \left(\mu_1 \mR_1 \mu _3+\eta _1 \mu _0\right)+s_0 \mu_0 \eta _1 \left(\mu_2 \mR_2 \eta _1-\mu_1 \mR_1 \eta _2\right)}{\left(\mu_2 \mR_2 \eta _1-\mu_1 \mR_1 \eta _2\right) \left(\eta _1 \mu _2-\eta _2 \mu _1\right)}}.
\eea

Until a point, the stability analysis is easy. DFE is stable iff $R_i \leq 1, i=1,2,3$.

\beR The formulas suggest introducing also
$$R_{13}=R_1(s_{13}),R_{23}=R_2(s_{23}),R_{31}=R_3(s_{13}),R_{32}=R_3(s_{23}),$$
where we recall that $s_{13}, s_{23}$ represent the boundary value of $s$ on $i_1=0$ and $i_2=0$, \resp.

\eeR
\sssec{Stability analysis for E1, E2}

Direct analysis implies that E1 is stable iff
\be{E1s}\bc R_1 \geq \max[1,R_2],\\ 
R_{13}=R_1(s_{13})>1 \Eq s_{13} >s_1,\ec\ee
with analog conditions holding for the stability of E2.
\beR
Recall that E13 exists iff $R_{13}>1$, and $R_{31}<1$,  thus stability of E1 does not prelude existence of E13. \eeR

To check \eqr{E1s}, note that the jacobian at $E1$ is $\Scale[0.82]{\left(
\begin{array}{cccc}
 -\mu _0 R_1 & -\mu _1 & -\frac{\mu _2 R_2}{R_1} & -\frac{\mu _3 R_3}{R_1} \\
 \mu _0 \left(R_1-1\right) & 0 & 0 & -\frac{\eta _1 \mu _0 \left(R_1-1\right) s_0}{\mu _1 R_1} \\
 0 & 0 & \mu _2 \left(\frac{R_2}{R_1}-1\right) & 0 \\
 0 & 0 & 0 & \frac{\mu _1 \mu _3 \left(R_3-R_1\right)+\eta _1 \mu _0 \left(R_1-1\right) s_0}{\mu _1 R_1} \\
\end{array}
\right)}$, and the numerator of the  \ch\ at $E1$  is:
$$\left(\mu _2 (R_1- R_2)+R_1 u\right) \left(\mu _0 \mu _1 (R_1-1)+\mu _0 R_1 u+u^2\right) \left(\mu _1 \mu _3 (R_1- R_3)+\eta _1 \mu _0 s_0(1- R_1)+\mu _1 R_1 u\right)$$

The linear terms are stable iff $$\Scale[.82]{R_1- R_2>0,\mu _1 \mu _3 (R_1- R_3)+\eta _1 \mu _0 s_0(1- R_1)>0 \Eq \mu _1 \mu _3 (R_1- R_3)+\eta _1 \mu _0 s_0<\eta _1 \mu _0 s_0R_1 \Eq \mR_1 s_{13}>1,}$$
and the quadratic term is stable iff $R_1-1>0$.

\beR Besides the  first  stability conditions for E1 in \eqr{E1s}, we have obtained a second one, which we are tempted to describe as ``E1 can invade E3", with the caveat that the invasions in situations with at least three siphons do  not seem to be rigorously defined -- see though \cite{Mitchell,Madec,kribs2023impact,belluccinimathematical} for some works in this direction.
\eeR
\sssec{Stability analysis for E3}
At E3, the \num\ of \ch\ is $\left(\mu _0 \mu _3 (R_3-1)+\mu _0 R_3 u+u^2\right)*$
$$
\left(\mu _1 \mu _3 (R_3- R_1)+\eta _1 \mu _0  s_0(R_3-1)+\mu _3 R_3 u\right) \left(\mu _2 \mu _3 (R_3- R_2)+\eta _1 \mu _0  s_0(R_3-1)+\mu _3 R_3 u\right).$$

Thus, stability holds iff
\bea
\bc \max[R_1,R_2,1]\leq R_3\\ R_{31}=\mR_3 s_{13}>1, R_{32}=\mR_3 s_{13}>1\ec\eea
holds,
   whose invasion interpretation is that besides $i_{12}$ being the dominant strain, it must also be able to invade strains $i_{1},i_{2}$.

\sssec{Stability analysis for E23}
At E13,E23 the direct \RH\ analysis becomes more cumbersome, because one of the factors of the \ch\ has degree 3,  but the  ``stability via LCP"
is not that hard.  Recall
$E23=\left\{s_{23}= \frac{s_0 \mu_0 \eta _2}
{\mu_2  \mu _3(\mR_2-\mR_3)+\eta _2 \mu _0},i_1= 0,
i_2=\frac{\mu _3}{\eta _2}(1-  \mR_3 s_{23}),i_3= \frac{\mu _2}{\eta _2}( \mR_2 s_{23}-1)\right\}$. To the positivity conditions of the above, we must add the negative flow condition near $i_1=0,$ which is
$ \mR_1 s<\fr{\eta _1}{\mu_1} i_3+1.$
A common reduce of these equations yields
\paragraph{Existence and Stability conditions for E23}
\begin{equation}
\bc \alpha_2\mu_3 > \alpha_3\mu_2 \Eq \mathcal{R}_2 > \mathcal{R}_3 \Eq R_2 > R_3\quad \text{and} \\
 \frac{\mu_2}{\alpha_2} < s_{23} < \frac{\mu_3}{\alpha_3} \Eq s_2 < s_{23}<s_3 \Eq  1< \mR_2 s_{23}=R_{23}, \mR_3 s_{23}=R_{32} <1\\
s_{23} < \fr{\eta_2\mu_1 - \eta_1\mu_2}{-\alpha_2\eta_1 + \alpha_1\eta_2}=\sE\ec
\end{equation}

The first three conditions ensure the existence. The fourth condition, ensuring stability, does not seem to be expressable  in terms of reproduction functions.
It is however elegantly expressed in terms of the $s$ values at E23 and EE.

The formulas above suggest:
\beO[partition via  resident variable  inequalities] \label{o:ins}
The parameter space may be divided via
 inequalities of the form $s_S\leq s_{S'}$  where $S \subset S'$  run over certain pairs of siphons (it may be possible to make this more precise using the theory of invasion graphs \cite{HofSch,Alm,Sch}.
\eeO

Finally, we show that this model does not \saty\ the CEP, \ie\ the existence conditions of EE are not that all the other fixed points exist and are unstable.
To simplify the existence conditions, assume first \wlo\ that
$$\eta_1/\mu_1 > \eta_2/\mu_2 \Lra \mR_2 <\mR_1 < \mR_2 \fr{\eta_1/\mu_1} { \eta_2/\mu_2}.$$ The existence conditions are still too long to reproduce here, so we further assume that all death rates equal $\mu$. Then we find a conjunction of four conditions
\bea &&1<R_2< R_1<1+ (R_2-1) \frac{\eta _1}{\eta _2}, R_3<\frac{\eta _1 R_2-\eta _2 R_1}{\eta _1-\eta _2},\mu >\sqrt{\frac{\eta _2 R_1}{\eta _1 R_2}}\\&&  \frac{\left(\eta _1-\eta _2\right) \mu  \left(R_1-R_3\right)}{\eta _1 \left(\eta _2+\eta _1 \left(R_2-1\right)-\eta _2 R_1\right)}<s_0<\frac{\left(\eta _1-\eta _2\right) \mu  \left(R_2-R_3\right)}{\eta _2 \left(\eta _1 (R_2-1)-\eta _2 \left(R_1-1\right)\right)},\eea
which leave several regions not covered by any of the previous cases.


\section{A Review of Persistence Theory for Positive Dynamical Systems}\lbl{s:per}

Let $\dot{x} = f(x)$ be a positive ODE defined on a compact, forward-invariant
set $E \subset \mathbb{R}_{+}^n$, with semiflow $\Phi_t$.

\begin{definition}[Persistence] \cite{TS,CraNazPan}.
The system is \textbf{persistent} if for every component $i$ and every trajectory starting in $\mathrm{int}(E)$,
\[
\liminf_{t \to \infty} x_i(t) > 0.
\]
\end{definition}

\begin{definition}[Uniform Persistence]\cite{TS,CraNazPan}.
The system is \textbf{uniformly persistent} if there exists $\varepsilon > 0$ such that
\[
\liminf_{t \to \infty} x_i(t) \ge \varepsilon
\]
uniformly for all trajectories starting in $\mathrm{int}(E)$.
\end{definition}

The fundamental problems of persistence, permanence, and extinction have been studied extensively, in particular in the context of generalized Lotka--Volterra systems
\cite{hutson1984average,butler1986persistence,hofbauer1989uniform,hutson1992permanence,
hofbauer1998evolutionary}, for mathematical ecology models \cite{thieme1992persistence,thieme2003mathematics,smith2011dynamical} and also for CRNs \cite{AdLS,JohnSie,PantGAS,CraNazPan,AndEncJoh,GopMilShiu}.

\subsection{General Persistence Theory}

\begin{theorem}[Butler--Waltman \cite{butler1986persistence}; Fonda \cite{fonda1988}; Hutson--Schmitt \cite{hutson1992permanence}; Smith--Thieme \cite{smith2011dynamical}]
\label{thm:repeller}
For a dissipative semiflow $\Phi_t$ on a compact invariant set $E$,
the system is uniformly persistent
if and only if the boundary $\partial E$ is a uniform repeller, i.e.,
there exists $\eta > 0$ such that
\[
\mathrm{dist}(\Phi_t(x), \partial E) \ge \eta
\]
for all $t$ sufficiently large and all $x \in \mathrm{int}(E)$.
\end{theorem}

\begin{remark}[On the insufficiency of linear instability for boundary fixed points]
The criterion above is not constructive. The simplest constructive condition—that all boundary equilibria are  unstable—does \emph{not} imply persistence -- see the May-Leonard example, subsection \ref{s:ML}.
This observation motivated the repeller theorem \ref{thm:repeller} of
Hutson--Schmitt~\cite{hutson1992permanence} and
Smith--Thieme~\cite{smith2011dynamical},
which assumes that the boundary repelling property must hold for the \emph{entire boundary set},
not just for fixed points.
\end{remark}

\subsection{Mass-Action CRNT Persistence Results}

In CRNT, the driving force behind persistence theory has been the global attractor conjecture 
for weakly-reversible systems, still unresolved; however, 
several  partial results are available \cite{AdLS,AndGAS,CraNazPan,PantGAS,CraGAS,GopMilShiu}.

\begin{theorem}[Angeli--De Leenheer--Sontag \cite{AdLS}, 2007]
\label{thm:angeli-sontag-persistence}
A mass-action chemical reaction network such that every siphon $Z$ of the network contains the support of at least one conservation law is persistent.
\end{theorem}

\begin{theorem}[Anderson \cite{AndGAS}, 2011]
\label{thm:anderson-persistence}
A weakly reversible, deficiency zero reaction network with mass-action kinetics \sats:
\begin{enumerate}
\item within each positive stoichiometric compatibility class, there exists a unique equilibrium;
\item this equilibrium is globally asymptotically stable;
\item the system is persistent on $\mathbb{R}^n_{++}$.
\end{enumerate}
\end{theorem}

\begin{theorem}[Pantea \cite{PantGAS}, 2012]
\label{thm:pantea-low-dimensional}
Every weakly reversible mass-action system of dimension $d \leq 3$ (where $d = \dim(S)$ is the dimension of the stoichiometric subspace) is persistent and permanent.
\end{theorem}

\beR
Pantea's proof uses:
\begin{enumerate}
\item Brouwer degree theory for $d = 2$
\item Poincaré--Bendixson theory ruling out boundary attractors
\item For $d = 3$: Analysis of possible boundary limit sets using $\omega$-limit set structure
\end{enumerate}
The result does \emph{not} extend to $d \geq 4$ where more complex boundary attractors can exist.
\eeR

\beR[Why standard CRN persistence results fail for epidemics]
Typical compartmental epidemic models satisfy:
\begin{enumerate}
\item \textbf{No weak reversibility}: Infection reactions $S + I \to 2I$ lack reverse $2I \to S + I$.
\item \textbf{Infection strains do not contain positive conservation laws}, due to possible deaths.
\item \textbf{Deficiency $\geq 1$}.

\end{enumerate}

\eeR

\subsection{Failure of Persistence due to Boundary Heteroclinic Cycles:May--Leonard system; MayLeonard.py, MayLeonard.nb}\lbl{s:ML}

Persistence may fail even when all boundary equilibria are  unstable.
One  mechanism illustrating this is an attracting \emph{boundary heteroclinic cycle}:
a sequence of boundary saddle equilibria connected cyclically,
each unstable toward the next, collectively forming a closed boundary orbit
that attracts all interior trajectories.
\beXa[Symmetric May--Leonard system]
\label{ex:mayleonard-ode}
A classical three-species competitive Lotka--Volterra equations
(\cite{may1975,hofbauer1993})
\begin{equation}\label{eq:mayleonard-ode}
\begin{aligned}
\dot{x}_1 &= x_1\big(1 - x_1 - \texttt{al}\,x_2 - \texttt{be}\,x_3\big),\\
\dot{x}_2 &= x_2\big(1 - \texttt{be}\,x_1 - x_2 - \texttt{al}\,x_3\big),\\
\dot{x}_3 &= x_3\big(1 - \texttt{al}\,x_1 - \texttt{be}\,x_2 - x_3\big).
\end{aligned}
\end{equation}
\sats\ that for $\texttt{al}<1<\texttt{be}$, the boundary equilibria
$E_1=(1,0,0)$, $E_2=(0,1,0)$, and $E_3=(0,0,1)$
are saddles connected in a heteroclinic cycle
$E_1\!\to\!E_2\!\to\!E_3\!\to\!E_1$
that attracts all interior trajectories, causing persistence failure.

\eeXa

The Figure \ref{fig:may-leonard-boundary} below  visualizes the May–Leonard heteroclinic cycle regime with parameters $al = 0.8$ and $be = 1.2$, demonstrating nonpersistence through three complementary perspectives:

\BEN

\im The left panel plot tracks $\log_{10}\min{x_1, x_2, x_3}$ over time essentially monitoring the smallest population abundance. The y-axis uses a logarithmic scale to capture orders of magnitude declines.
Starting from multiple initial conditions (shown as different colored trajectories), all paths exhibit monotone decay toward $-\infty$ on the log scale. This means the minimum population size continually decreases, approaching zero asymptotically. The smooth downward trend confirms that no trajectory stabilizes at a positive interior equilibrium; instead, populations drift inexorably toward boundary faces where at least one species becomes vanishingly rare.

\im For the middle plot, it shows a single representative trajectory projected onto the 2-simplex (the triangular constraint $x_1 + x_2 + x_3 = 1$), color-coded by distance to the nearest boundary face.
\BEN
\im Yellow bright regions: trajectory far from boundaries (interior of simplex),
\im Blue dark regions: trajectory near boundaries (edges or vertices).
\EEN
Darker shading near the vertices indicates prolonged visits to states where one species dominates temporarily while the other two remain suppressed. The trajectory never settles but instead circulates endlessly around the boundary loop, never stabilizing in the interior.

\im For the right pannel, it plots the time evolution of all three species over the interval $[0, 2000]$ time units, capturing long-term oscillatory behavior.
The vertical axis shows population densities $x_1, x_2, x_3$, and horizontal lines separate successive epochs where different species dominate.
\EEN

This multi-faceted visualization clarifies why classical stability theory focused on interior equilibria and local linearizations fails to predict dynamics in systems with boundary attractors. The May-Leonard example remains a canonical warning that ecological coexistence requires more than merely unstable single species states; the connectivity and geometry of boundary equilibria matter profoundly.

\begin{figure}[H]
\centering
\begin{minipage}{0.32\linewidth}
  \includegraphics[width=\linewidth]{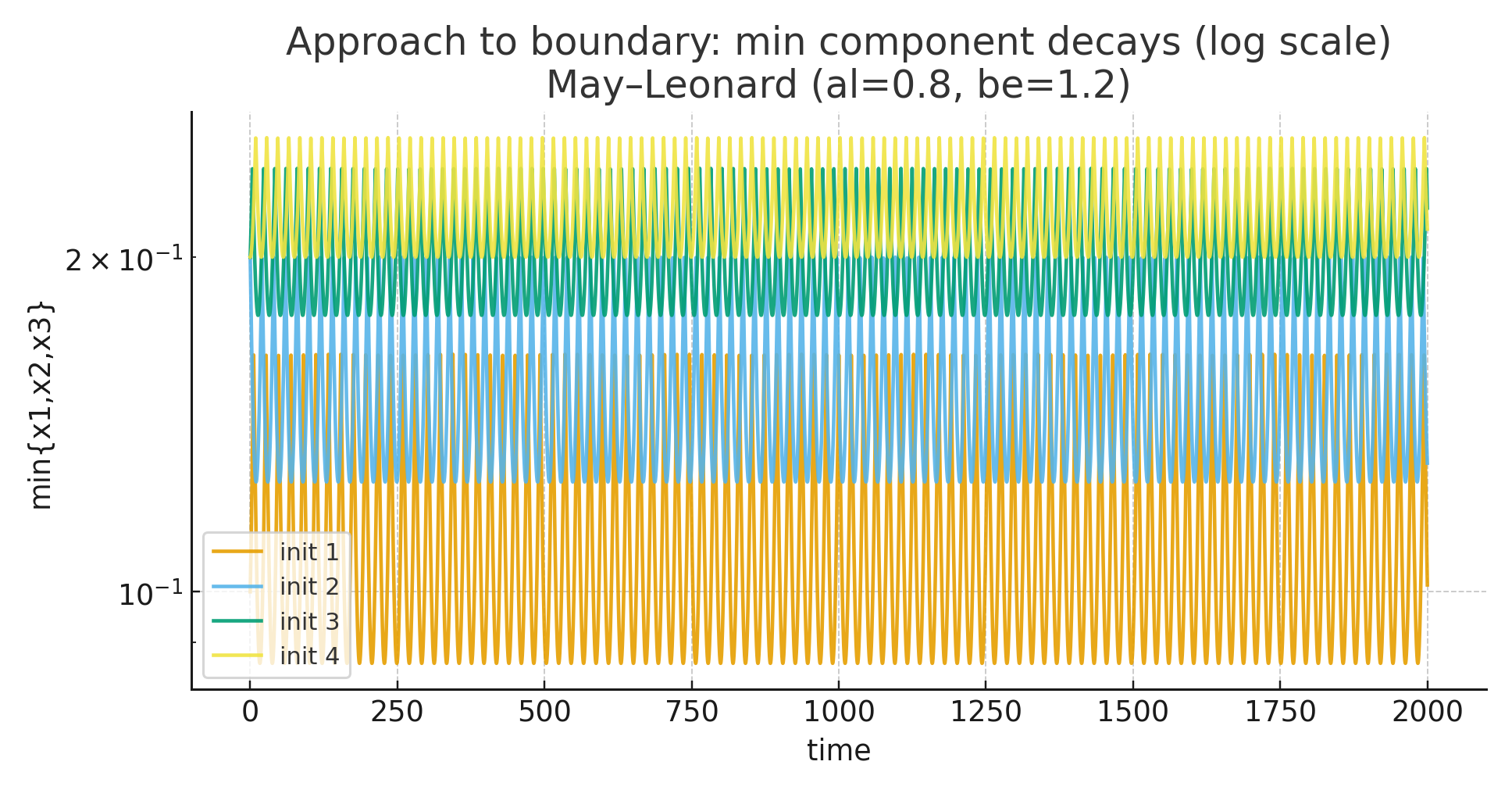}
  \caption*{\small $\log_{10}\!\min\{x_1,x_2,x_3\}$ vs.\ time for several initial conditions.
  Monotone decay evidences approach to the boundary.}
\end{minipage}\hfill
\begin{minipage}{0.32\linewidth}
  \includegraphics[width=\linewidth]{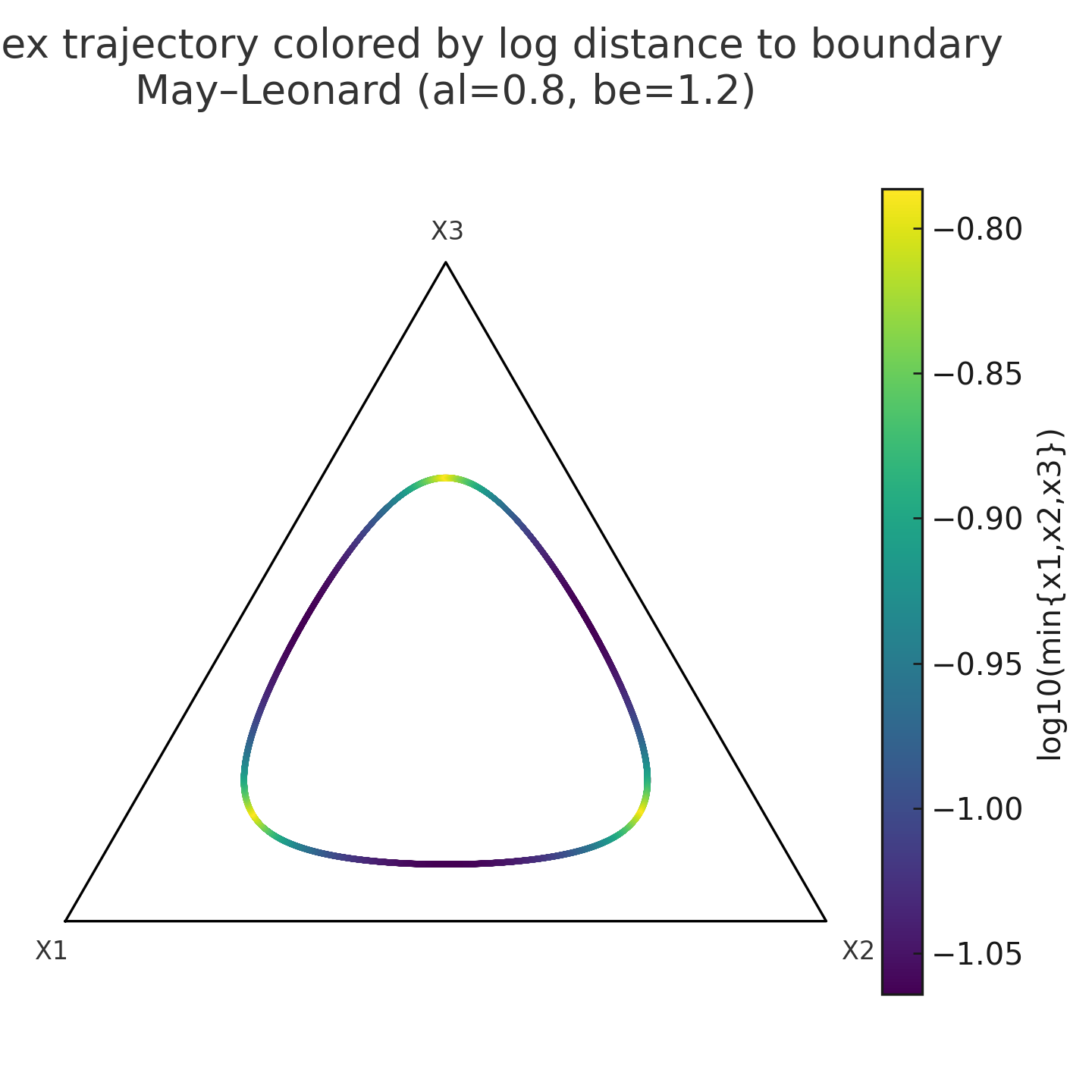}
  \caption*{\small Simplex trajectory colored by $\log_{10}\!\min\{x_1,x_2,x_3\}$.
  Darker colors indicate closer proximity to boundary faces.}
\end{minipage}\hfill
\begin{minipage}{0.32\linewidth}
  \includegraphics[width=\linewidth]{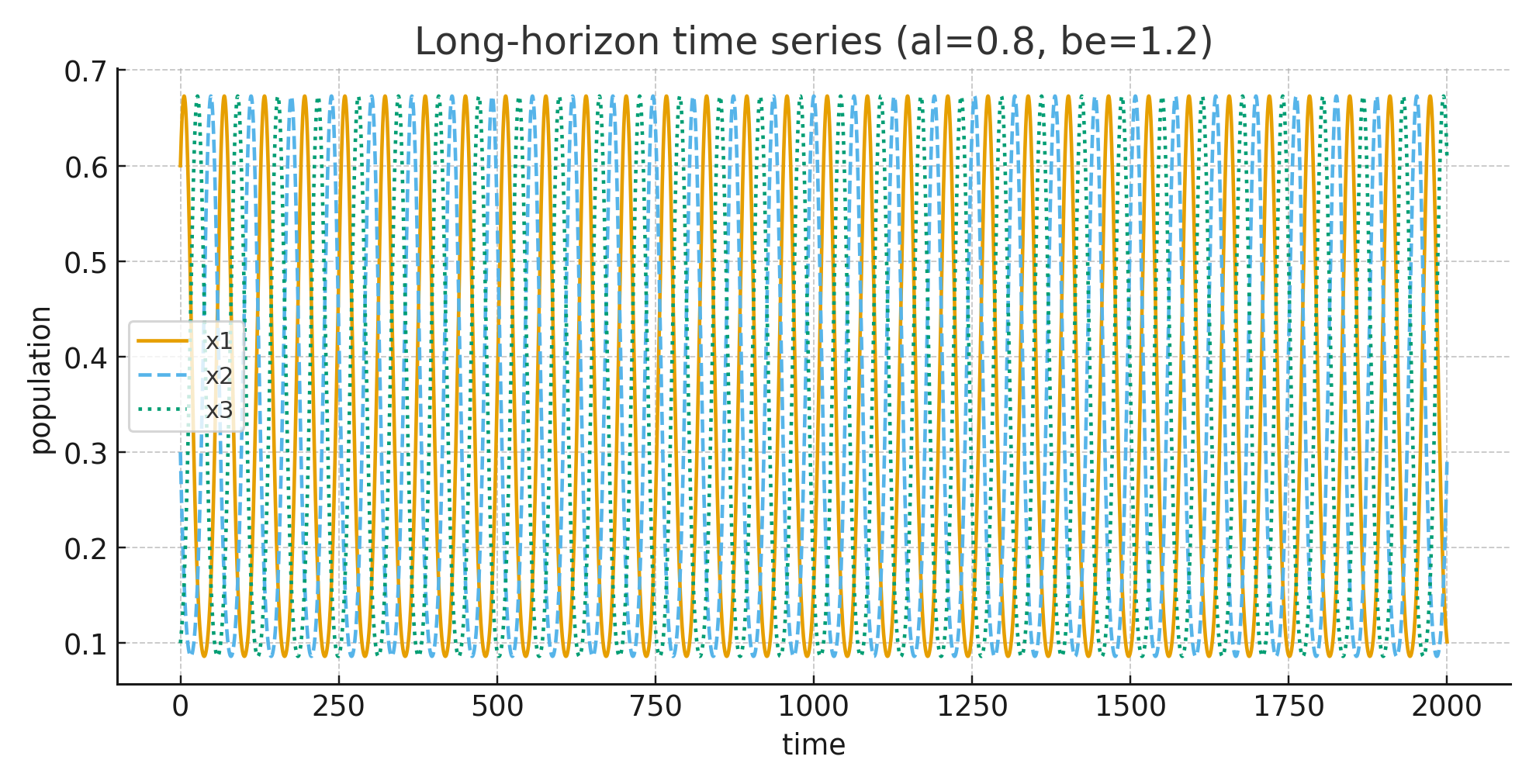}
  \caption*{\small Long-horizon time series showing cycling dominance and longer boundary visits.}
\end{minipage}
\caption{May–Leonard heteroclinic-cycle regime with \texttt{al}$=0.8$, \texttt{be}$=1.2$.
Left: $\log_{10}\!\min\{x_1,x_2,x_3\}$ clearly $\to -\infty$, hence nonpersistence.
Middle: simplex projection with distance-to-boundary coloring.
Right: time series on $[0,2000]$ showing the characteristic increase of dwell times near the boundary saddles. These are obtained in MayLeonard.nb}
\label{fig:may-leonard-boundary}
\end{figure}

For checking this with EpidCRN, we use the
mass--action  representation:
\begin{verbatim}
RN = {
  "X1" -> 2*"X1",
  "X1" + "X1" -> "X1",
  "X1" + "X2" -> "X2",
  "X1" + "X3" -> "X3",

  "X2" -> 2*"X2",
  "X2" + "X2" -> "X2",
  "X2" + "X3" -> "X3",
  "X2" + "X1" -> "X1",

  "X3" -> 2*"X3",
  "X3" + "X3" -> "X3",
  "X3" + "X1" -> "X1",
  "X3" + "X2" -> "X2"
};
rts = {
  x1, x1*x1, al*x1*x2, be*x1*x3,
  x2, x2*x2, al*x2*x3, be*x2*x1,
  x3, x3*x3, al*x3*x1, be*x3*x2
};
\end{verbatim}
\subsection{Invasion graphs}

In this subsection we revisit the invasion graphs theory which arose from the works of
\cite{Schreiber,Schreiber06,PatelSchreiber,Mitchell,Madec,HofSch,LeGjiniMadec2023,Schreiber25}
using the notions of siphons and reproduction functions.

Throughout, $X\subset\R^{m+n}_+$ is the phase space of an ME system with
infection set $TS=\{1,\dots,n\}$.

\beD[Proper siphons and communities]\label{d:prS}
Let $S\subset TS$ be a siphon and let
\[
F_S=\{x\in X:\ x_j=0 \ \text{for all } j\in S\}
\]
denote the associated positively invariant face.

A siphon $S$ is called {\em proper} if the restriction of the flow to $F_S$
admits a compact invariant set $M_S\subset F_S$ and $M_S$ contains a point with all coordinates in $S^c$ strictly positive.

The complement $S^c=TS\setminus S$ is called the {\em community}.
\eeD

\beD[Invasion rates]\label{d:inR}
Let $S$ be a proper siphon and $M_S$ its compact invariant set.
For $j\in S$ the $j$--equation linearizes along $\phi_t(x)$, $x\in M_S$, to
\[
z'_j = a_j(t)\,z_j,
\qquad
a_j(t)=\partial_{x_j}f_j(\phi_t(x)).
\]
The (upper) Lyapunov exponent
\[
\lambda_j(M_S)=
\limsup_{t\to\infty}\frac{1}{t}\int_0^t a_j(\tau)\,d\tau
\]
is the {\em invasion rate} of coordinate $j$ into the community $S^c$.
We say $j$ {\em invades} $S^c$ whenever $\lambda_j(M_S)>0$.
\eeD

\begin{remark}[Reduction for ME models with resident equilibria]
In all ME models considered here, $M_S$ is a single equilibrium
$E_{S^c}$ with susceptible level $s_S$.  Then $a_j(t)\equiv a_j$, and
\[
\lambda_j(M_S)=\partial_{x_j}f_j(E_{S^c}),
\]
so invasion rate reduces to a single Jacobian evaluation.

Moreover, using the $F{-}V$ splitting of infection terms,
\[
\lambda_j(M_S)>0 \Longleftrightarrow R_j(s_S)>1
\]
where $R_j$ is the reproduction function of the $j$--th NGM block.
\end{remark}

\beD[Invasion numbers for ME models]\label{d:inv3}
Let $T_j$ be an NGM block, with reproduction function $R_j(\cdot)$.
For a proper siphon $S$, let $E_{S^c}=(x_S=0,\ y_{S^c})$ be the resident
equilibrium on $F_S$, and let $s_S$ denote its susceptible coordinate.

The {\em invasion number} of block $T_j$ into the community $S^c$ is
\[
R_j^{S^c}=R_j(s_S).
\]
The {\em reproduction number} of $T_j$ is
\[
R_j^0=R_j(s_0)=R_j(E_0).
\]
\eeD

\begin{remark}
If the minimal siphons partition $TS$ and determine the NGM blocks, the
preceding definition reduces to the previous definition \ref{d:reF}.
\end{remark}

\beXa[Invasion graph of the GK coinfection model]
\label{ex:inGK}
We consider the GK model with blocks
\[
T_1=\{i_1\},\qquad T_2=\{i_2\},\qquad T_3=\{i_{12}\},
\]
and reproduction functions
\[
R_j(s)=\mathcal{R}_j s,\qquad \mathcal{R}_j=\frac{\alpha_j}{\mu_j}.
\]

For $\beta>0$ and $\gamma>0$ the minimal siphons are
\[
\Sigma_1=\{i_2,i_{12}\},\qquad \Sigma_2=\{i_1,i_{12}\},
\quad TS=\{i_1,i_2,i_{12}\},
\]
and the proper siphons are
\[
S_*=\emptyset,\qquad S_1=\Sigma_1,\qquad S_2=\Sigma_2,\qquad S_0=TS.
\]

Each $S$ determines a unique resident equilibrium $E_{S^c}$ with
susceptible value $s_S$.  Table~\ref{tab:sGK} lists all these values.
Invasion inequalities are of the form
\[
R_j(s_S)>1,\qquad j\in S.
\]

\begin{table}[h]
\centering
\caption{Resident susceptible values $s_S$ for proper siphons in the GK model.}
\label{tab:sGK}
\begin{tabular}{c|c|c}
Proper siphon $S$ & Community $S^c$ & Resident susceptible value $s_S$\\
\hline
$S_0=\{1,2,3\}$ & $\emptyset$ & $s_0$ (DFE value)\\
$S_1=\{2,3\}$ & $\{1\}$ & $s_1=\dfrac{1}{\mathcal{R}_1}$\\
$S_2=\{1,3\}$ & $\{2\}$ & $s_2=\dfrac{1}{\mathcal{R}_2}$\\
$S_*=\emptyset$ & $\{1,2,3\}$ & $s_*=s_3=\dfrac{1}{\mathcal{R}_3}$\\
\end{tabular}
\end{table}
Note that the indices $\pr{1,2,3}$ correspond to $\pr{i_1,i_2,i_{12}}$.

To handle mixed equilibria on $i_1=0$ or $i_2=0$ we also introduce
\[
s_{13}=\pi_s(E_{13}),\qquad s_{23}=\pi_s(E_{23}),
\]
given explicitly in the main text.
All invasion inequalities $R_j(s_S)>1$ are listed in
Table~\ref{tab:invGK}.

\begin{table}[h]
\centering
\caption{Invasion inequalities $R_j(s_S)>1$ for GK.
Each inequality corresponds to a directed edge $S\to S'$ in the invasion graph.}
\label{tab:invGK}
\begin{tabular}{c|c|c}
Invaded siphon $S$ & Invading block $T_j$ & Inequality $R_j(s_S)>1$\\
\hline
$S_0$ & $T_1$ & $R_1(s_0)>1$\\
$S_0$ & $T_2$ & $R_2(s_0)>1$\\
$S_0$ & $T_3$ & $R_3(s_0)>1$\\
\hline
$S_1$ & $T_2$ & $R_2(s_1)>1$\\
$S_2$ & $T_1$ & $R_1(s_2)>1$\\
$S_1$ & $T_3$ & $R_3(s_{13})>1$\\
$S_2$ & $T_3$ & $R_3(s_{23})>1$\\
\hline
$S_*$ & $T_1$ & $R_1(s_*)>1$\\
$S_*$ & $T_2$ & $R_2(s_*)>1$\\
$S_*$ & $T_3$ & $R_3(s_*)>1$
\end{tabular}
\end{table}

The {\em invasion graph} is the directed graph whose vertices are the
proper siphons
\[
S_*,\ S_1,\ S_2,\ S_0
\]
and whose edges are exactly those pairs $S\to S'$ for which the block
removed in passing from $S$ to $S'$ satisfies $R_j(s_S)>1$.
\eeXa


\section{Appendix A}\lbl{s:ap}

\subsection{Glossary of  the acronyms}

\begin{center}
\begin{tabularx}{\linewidth}{|l|X|X|}
\hline
\textbf{Acronymes} &
\textbf{Explanation} &
\textbf{References} \\
\hline
MS/AMSD & Minimal Siphon/ Acyclic Minimal Siphon Decomposition & \cite{ShiuStu, you2021}
\\\hline
IGMS & interaction graph on minimal siphons &  \cite{ShiuStu,banaji2009graph}\\ \hline
SDAS & Seed-Dependent Autocatalytic Systems & \cite{Peng2021SDAS,Peng2022HierarchicalSDAS,Peng2023ChemicalEcosystem,
 Blok23} \\\hline
LCP/CEP &  linear complementarity principle/competitive exclusion principle & \cite{thieme89, camlibel2007} \\\hline
DSR & Directed Species-Reaction graph &  \cite{banaji2009graph}\\\hline
NGM & Next-Generation Matrix & \cNGM\ \\\hline
\end{tabularx}
\end{center}

\subsection{A useful background result for positive ODEs : the Regular Splitting theorem for Metzler matrices}\label{s:reg}

The NGM method  consists of two steps: \BEN \im   the reduction of the Jacobian to the \jin, and, assuming this is Metzler,
 \im  studying its stability by finding a \emph{regular matrix splitting}.
 \EEN
 \beD[regular splitting] \label{d:RS} \cite{Varga,Plemmons}. Given a Metzler (quasi-positive) matrix $M$, a decomposition $M = F -V$ is a regular splitting if:
\begin{itemize}
    \item $F \geq 0$
    \item $V^{-1} \geq 0$
\end{itemize}

Equivalently, $A=-V$ must be a \MH\ /Markovian semi-group generating  matrix.\eeD

\beR Note  that
a \dec\ like in Definition \ref{d:RS} needs neither exist, nor be unique.

In quadratic ME examples, a regular splitting is sometimes obtained by  separating the quadratic part
$F$ of $J_x$ from the linear part $V$.
\eeR
 Assuming a regular splitting exists,
the spectral radius of the \emph{next generation matrix}:
\[
K := F(-A)^{-1}
\]
provides a representation of   the stability domain of $M$:
\beP \label{p:Var} (Varga, \cite[Thm 3.13]{Varga}, \cite{Berman,Fall, thieme2009spectral})

 Let $M$ be a real Metzler matrix  admitting regular splitting(s).
 Then
 $$\bc \rho ( F V^{- 1}) < 1\Lra M &\text{is \MH\ }\\
 \rho ( F V^{- 1}) > 1\Lra M &\text{is unstable }\end{cases},$$ where $M =  F-V$ is  any regular splitting of  $M $.
\eeP
This result originates in Varga's work on iterative methods.

\beD [next generation matrix] \cite{Diek10,Van}. For   a real Metzler matrix $M$ admitting a regular splitting
$M=F-V,$ the matrices $F V^{- 1}, V^{- 1} F$ (which have equal spectrum) will both be called next generation matrices. \eeD

Below, we will use the first form.

\ssec{Informal description of the  next generation matrix\ method}\label{s:NGM}

Once the DFE is identified, its  stability
 threshold   may be determined via the NGM method, which analyzes the DFE by projecting the Jacobian onto infection variables. Let $M = \partial f / \partial x$ at the DFE. If $M$ is Metzler and admits a regular splitting $M = F - V$, the DFE is locally stable iff:
\[
R_0 =R_0(F):= \rho(K): = \rho(F V^{-1}) < 1
\]

 We  give now an informal description of the NGM method (see \cNGM\ for precise statements):
 \BEN \im One restricts to the "disease \equ s", \ie\ the equations involving the $0$ coordinates of the DFE. This  may be justified by
 \BEN
 \im verifying the conditions of \cNGM, or, in particular cases, by\im
 factoring the \ch\ of the Jacobian matrix,  and by eliminating trivial factors (which correspond always to the non-disease \var s).

 \EEN

 \im  Assuming the \jin\ $J_x$ is a Metzler matrix,
 and that it admits a regular splitting \be{rs}J_x=F +A:=F-V,\ee the  study of  $J_x$ may be replaced by that of one of the "next generation matrices"
  \be{K} K_d= V^{-1} F, K= FV^{-1} =
 V K_d V^{-1}\ee
 which have probabilistic interpretations, and may also be obtained by an algebraic trick:
 \be{Kl}J_x=F -V=V\pp{V^{-1}F  -I}=  V\pp{K_d-I}=\pp{K-I} V.\ee

   The final result is the DFE stability threshold
 \be{R0fr}R_0=R_0(F):=\rho(K)=\rho(K_d)<1\ee
 (a corollary of Varga's regular splitting Lemma \cite[Thm 3.13]{Varga}, \cite{Berman,Fall, thieme2009spectral}).

 \EEN

We will call the inequality \eqr{R0fr} the first law of \ME; it states that the (interior of the) stability domain of the DFE has a representation \eqr{R0fr}, for any admissible matrix $F$.

 \beR Note that we introduced above two NGMs, $K_d$, and $K$, related by a similarity transformation, and
 having thus the same  \ch, $R_0$, and   stability properties. 

The comparison of the two matrices
 $K,K_d$   was addressed in \cite{Diek10}, who  claimed  that $K_d$ leads to simpler computations  than $K$, but we have encountered mainly the opposite situation, and prefer to use $K$. Let us  add
here an empirical observation, that the
computation of $R_0$ from $K$ or $K_d$ is typically easier than via the reduced Jacobian $J_x$, in the sense that their characteristic polynomials factor more.

Further research into   the respective advantages of the two NGMs seems necessary. 

\eeR

\beD[minimal $R_0$ degree] When $R_0$ is algebraic, the degree of the smallest irreducible polynomial \satd\ by  $R_0(F)$
for some \adm\ \dec\ will be called the degree of the  \dec, and the
minimal  degree over all \adm\ \dec s will be called the minimal $R_0$ degree.\eeD

\ssec{A heuristic NGM approach}
We believe it  may be
 useful to complement  the rigorous NGM \cite{Van} with  the following
``heuristic" NGM approach (which is actually what many papers seem to use in practice):
\begin{enumerate}
    \item For a positive ODE,  compute the union of minimal siphons $\x_{E}$, and verify algebraically the existence of a  boundary fixed point\ $E$ on the facet determined by $\x_{E}$.
        \im Compute the resident Jacobian block $J_y(0,\y_E)= \fr{\partial f }{ \partial \y}(0,\y_E)$ and verify its stability.
    \item If it is stable, compute also the invasion Jacobian block $$J_x(0,\y),$$  without plugging the values of the resident \var s (this is useful for symbolic manipulations, and has also theoretical implications which require further investigation), and verify whether it is Metzler.

    \item If it is, find a decomposition $J_x = F - V$, where:
        \begin{itemize}
            \item $F$: gain terms (nonlinear, positive)
            \item $V$: loss/removal terms verifies that $V^{-1}$ has only positive elements.
        \end{itemize}
        \end{enumerate}

        If all the steps were successful, we are dealing with an ME model,
        and we may further compute the matrix $K = FV^{-1}$.
        \beR Sometimes,  eliminating trivial factors from the characteristic polynomial of $K$ further   eases  the symbolic analysis, effectively replacing the model by a simpler one. \eeR

         In the opposite case, we would have identified the reason why our  positive ODE  is not a ME model.

\beO [conditions for the $R_0$ dichotomy to hold, in terms of the partition  Jacobians] Provide conditions for the $R_0$ dichotomy to hold for ME models, in terms of the partition of the Jacobian in $J_x,J_y,J_{xy}, J_{yx}.$ \eeO

\section{Appendix B: Further pointers to the literature and ideas}

\subsection{Relation between minimal self-replicating siphons and autocatalytic cores}\lbl{s:rel}

\subsubsection{Background}

In the framework of persistence theory,
Deshpande and Gopalkrishnan~\cite{DesGop}
introduced the notion of \emph{self-replicable siphons}.
Let \(S\) denote the stoichiometric matrix of a chemical reaction network (CRN),
and let \(A\) be the restriction of \(S\) to a subset of species
forming a minimal critical siphon.
A siphon \(T\) is said to be \emph{self-replicable}
if there exists a nonnegative flux vector \(v>0\)
such that \(A v > 0\).
Theorem~5.3.3 of~\cite{DesGop} shows that
every minimal critical siphon is either \emph{drainable}
(\(A v < 0\) for some \(v>0\))
or \emph{self-replicable},
and that CRNs without drainable siphons are persistent.

Several years later, Blokhuis \emph{et~al.}~\cite{BlokLacNghe}
introduced the concept of \emph{autocatalytic cores} as
the minimal stoichiometric submatrices \(A\) satisfying
\(A v > 0\) for some \(v>0\),
subject to additional structural constraints:
each reaction must consume and produce species within the same
subnetwork and involve exactly one reactant per column.
These authors proved that
minimal autocatalytic matrices are necessarily
\emph{square, invertible, and Metzler},
thus linking stoichiometric minimality to algebraic irreducibility.

More recently,
Vassena and Stadler~\cite{VasSta}
formalized this connection in dynamical terms:
an \emph{autocatalytic core} are precisely the subset of
\emph{child-selection matrices} that are both
Metzler and Hurwitz-unstable,
which is equivalent to the existence of a positive vector \(v>0\) for which
\(A v > 0\).

\subsubsection{Definitions}
Let $\mathcal{G}=(\mathcal{S},\mathcal{R})$ be a chemical reaction network with stoichiometric matrix $S$.

\beD[Self-replicating siphon] \cite{DesGop}. 
A siphon $T\subseteq\mathcal{S}$ is \emph{self-replicating} if there exists a flux
$v\ge0$, $v\ne0$, such that $(S v)_i>0$ for all $i\in T$ and $(S v)_i=0$ for all $i\notin T$.
\eeD

The existence of such a flux ensures that, once $T$ is nonempty, the reactions of the
network can increase all species in~$T$ without aid from outside $T$.

\beD[Autocatalytic Siphon]\cite{DesGop,unterberger2022}.

Let $G = (S, R)$ be a chemical reaction network with stoichiometric matrix $S$. 
A siphon $T \subseteq S$ is called \emph{autocatalytic} if there exists a nonzero flux vector $v \ge 0$ such that:

\begin{enumerate}
    \item All species in $T$ increase under the flux:
    \[
    (Sv)_i > 0 \quad \text{for all } i \in T;
    \]
    \item Every species in $T$ participates as a reactant in at least one reaction contributing to $v$, meaning it helps catalyze its own production;
    \item The flux is supported entirely within $T$:
    \[
    (Sv)_i = 0 \quad \text{for all } i \notin T.
    \]
\end{enumerate}

In other words, once $T$ is nonempty, the species in $T$ can increase themselves through internal reactions, without relying on species outside $T$.
\eeD

\beR
\BEN
\im Every autocatalytic siphon is a self-replicating siphon, but not every self-replicating siphon is autocatalytic. 
\im Autocatalytic siphons capture the idea of \emph{internal catalytic feedback} within a siphon, which makes them key structures for persistence and potential growth in chemical reaction networks.
\EEN
\eeR

\beT[Deshpande--Gopalkrishnan]
\label{thm:DG}
\cite{DesGop}:Every minimal critical siphon of a weakly reversible CRN is either
\emph{drainable} or \emph{self-replicating}.
\eeT

\beT[Blokhuis, Lacoste, Nghe]
\label{thm:Blokhuis}
\cite{BlokLacNghe}: For mass-action kinetics, every minimal self-replicating siphon contains
a minimal \emph{autocatalytic core}, i.e.\ a minimal square submatrix
$A$ of~$S$ admitting a positive flux vector $v>0$ with $Av>0$.
Conversely, every autocatalytic core defines a self-replicating siphon
equal to its species support.
\eeT

\beT[Vassena--Stadler]
\cite{VasSta}: Minimal autocatalytic cores are precisely the minimal
Metzler–Hurwitz-unstable submatrices of the Jacobian at boundary equilibria.
Hence each minimal self-replicating siphon corresponds to a boundary-unstable
face of the positive orthant.
\eeT

\begin{remark}
Theorems~\ref{thm:DG}–\ref{thm:Blokhuis} establish a bijective correspondence
between minimal self-replicating siphons and autocatalytic cores:
the structural condition of self-replication (existence of a positive
stoichiometric flux restricted to~$T$) coincides with the dynamical
condition of local instability (positive Perron–Frobenius eigenvalue of a
Metzler block).
\end{remark}

\subsubsection{Comparative Summary}
Although the discussion above suggest a relation between
\emph{self-replicable siphons} and \emph{autocatalytic cores},
\textbf{precise results} in this direction are more read between the lines, than stated explicitly in the literature.
Conceptually, Deshpande’s self-replicable siphons describe
\emph{combinatorial minimal positive supports}
in the stoichiometric geometry,
while Blokhuis cores represent
\emph{minimal autonomous subnetworks}
that satisfy the same positivity condition together with
stoichiometric closure and invertibility.
In this sense, autocatalytic cores can be interpreted
as \emph{square, autonomous realizations}
of minimal self-replicable siphons,
but the equivalence remains heuristic.

\begin{center}
\begin{tabularx}{\linewidth}{X X X X}
\hline
\textbf{Concept} &
\textbf{Defining Condition} &
\textbf{Structural Constraints} &
\textbf{Source} \\
\hline
Self-replicable siphon &
$\exists\,v>0:\, A v > 0$ &
Minimal critical subset of species (not necessarily square) &
Desh. \& Gopal. (2014) \\[4pt]
Autocatalytic core &
$\exists\,v>0:\, A v > 0$ and $A$ square, invertible, Metzler &
Each species unique reactant; autonomous subnetwork &
Blokhuis \emph{et~al.} (2020) \\[4pt]
Unstable core (child selection) &
$A$ Metzler and Hurwitz-unstable $\Leftrightarrow \exists\,v>0:\,A v > 0$ &
Dynamical realization of an autocatalytic core &
Vassena \& Stadler (2024) \\
\hline
\end{tabularx}
\end{center}

\subsubsection{Discussion}

The absence of a formal equivalence theorem can be attributed
to methodological differences:
self-replicable siphons are defined combinatorially
in terms of reachable states and persistence,
while cores are algebraic objects derived from the
stoichiometric matrix structure.
Nonetheless, both capture the same necessary positivity condition
that guarantees self-maintenance and amplification within a subnetwork.


\subsection{Tree IGMS vs.  multi-tier Seed-Dependent Autocatalytic Systems (SDAS)}\lbl{s:SDAS}
\subsubsection{Definitions}

Let $\mathcal S$ be the set of species, $R$ the set of reactions, and $S\in\mathbb{R}^{|\mathcal S|\times |R|}$ the stoichiometric matrix.
\beD[Stoichiometric autocatalysis and core] \cite{DesGop}. 
For $U\subseteq\mathcal S$ and $R_U\subseteq R$, the restricted matrix $A=S_{U,R_U}$ is \emph{stoichiometrically autocatalytic} if there exists $v\ge0$, $v\ne0$, such that $Av>0$.
A minimal square such submatrix is an \emph{autocatalytic core}.
\eeD

\beD[Strict catalysis]\cite{DesGop}.  
A motif (flux vector $v\ge0$) is \emph{catalytic} if $Av>0$.
\eeD

\beD[formal autocatalysis]\cite{andersen2021defining}.
A subnetwork $\mathcal{G}' = (\mathcal{S}',\mathcal{R}')$ of CRN $\mathcal{G}= (\mathcal{S},\mathcal{R})$ is defined to be \textit{formally autocatalytic} in the subset $\mathcal{M}$ of $\mathcal{S}'$ if (see \cite{andersen2021defining}):
\begin{enumerate}
    \item There exists a positive real flow ($\gg 0$) on $\mathcal{G}'$ such that the resulting composite reaction is of the form
    \begin{align*}
        \text{(F)} + \textbf{m} \mathcal{M} \to \textbf{n} \mathcal{M} + \text{(W)},
    \end{align*}
    where $\textbf{0} \ll \textbf{m} \ll \textbf{n}$. Here $\textbf{m}$ and $\textbf{n}$ are the stoichiometries of the set $\mathcal{M}$ in the input and output of the composite reaction, respectively, and $\textbf{o} \mathcal{M} = \sum_i \textbf{o}_i \mathcal{M}_i$.
\end{enumerate}
\eeD

\beD[exclusive autocatalysis]\cite{andersen2021defining}.
\label{def:EA_GT}
A formally autocatalytic subnetwork $\mathcal{G}' = (\mathcal{S}',\mathcal{R}')$ of CRN $\mathcal{G}= (\mathcal{S},\mathcal{R})$ is defined to be \textit{exclusively autocatalytic} in the subset $\mathcal{M}$ (see \cite{DesGop,andersen2021defining}) if:
\begin{itemize}
    \item Every reaction in $\mathcal{R}'$ consumes at least one species from the set $\mathcal{M}$. This ensures that the flow is inadmissible, or there is no flow, if the population of any species in the set $\mathcal{M}$ is zero.
\end{itemize}
\eeD

\beD[Seed-dependent autocatalytic system (SDAS)]
\cite{Peng2022HierarchicalSDAS}. 
Let $\mathcal F\subseteq\mathcal S$ denote the externally supplied species.
A subnetwork $(\mathcal S',R')$ is an \emph{SDAS} if there exists a nonempty seed set $\Sigma\subseteq\mathcal S'\setminus\mathcal F$ such that
(i) $\mathcal F\cup\Sigma$ generates all reactants of $R'$ under $R'$, and
(ii) there exists $v\ge0$, $v\ne0$, with $S_{\mathcal S',R'}v>0$.
\eeD

\begin{proposition}[Tree IGMS is not sufficient for SDAS]
Acyclicity of IGMS is necessary for hierarchical seeding but does not imply the existence of a positive flux $v$ with $Av>0$.  Hence a tree IGMS does not guarantee an SDAS.
\end{proposition}

\begin{proposition}[Cyclic IGMS can admit SDAS]
Mutual seeding between two siphons may yield a cycle in IGMS while a joint autocatalytic core on their union satisfies $Av>0$ and constitutes an SDAS.
\end{proposition}

\subsubsection{Examples}

The following examples illustrate the correspondence between the
definition of \emph{exclusive autocatalysis} in a subset
$\mathcal M\subseteq\mathcal S$ (Definition~\ref{def:EA_GT}) and
\emph{stoichiometric autocatalysis} defined by the existence of a
nonnegative flux vector $v$ such that $S_{\mathcal M,R'}v>0$.
Condition~(i) of Definition~\ref{def:EA_GT} ensures that all reactions
in the subnetwork consume at least one member of~$\mathcal M$, so any
feasible flux must vanish if any component of $\mathcal M$ is zero.
Stoichiometric autocatalysis adds the quantitative requirement that a
strictly positive flux $v$ exists, implying net production of all
species in~$\mathcal M$.

In the following examples, $\mathcal M$ is indicated for each
subnetwork, and the feasibility of $Av>0$ determines whether it is
stoichiometrically autocatalytic.

\begin{example}[Tree IGMS without stoichiometric autocatalysis]
\begin{verbatim}
RN = {
  0 -> "S",
  "S" -> 0,
  "S" + "I1" -> "I1",
  "I1" + "I2" -> "I2",
  "I1" -> 0,
  "I2" -> 0
};
rts = {
  la*s,
  mu*s,
  be*s*i1,
  ro*i1*i2,
  mu1*i1,
  mu2*i2
};
\end{verbatim}
Here $\mathcal M=\{\text{"I1"},\text{"I2"}\}$ is exclusive: every reaction
consumes at least one member of~$\mathcal M$.
However, the restricted stoichiometric matrix $A_{\mathcal M}$ has no
nonzero $v\ge0$ with $A_{\mathcal M}v>0$.
The network is therefore exclusively but not stoichiometrically
autocatalytic, yielding a tree IGMS without an SDAS.
\end{example}

\begin{example}[Cyclic IGMS with stoichiometric autocatalysis]
\begin{verbatim}
RN = {
  0 -> "S",
  "S" -> 0,
  "S" + "I1" -> 2*"I1",
  "I1" + "I2" -> 2*"I2",
  "I2" + "I1" -> 2*"I1",
  "I1" -> 0,
  "I2" -> 0
};
rts = {
  la*s,
  mu*s,
  be*s*i1,
  ro12*i1*i2,
  ro21*i2*i1,
  mu1*i1,
  mu2*i2
};
\end{verbatim}
Each reaction consumes either $I_1$ or $I_2$, hence the subsystem
is exclusive in $\mathcal M=\{I_1,I_2\}$.
The restricted matrix $A_{\mathcal M}$ admits a vector $v>0$
with $A_{\mathcal M}v>0$, so the network is both exclusive and
stoichiometrically autocatalytic.
The corresponding IGMS contains a two-cycle
$T_1\!\leftrightarrow\!T_2$, demonstrating that cyclic IGMS
can still realize an SDAS.
\end{example}

\begin{example}[Three-tier SDAS with strict catalysis]
\begin{verbatim}
RN = {
  0 -> "S",
  "S" -> 0,
  "S" + "I1" -> 2*"I1",
  "I1" -> 0,
  "I1" + "I2" -> 2*"I2",
  "I2" -> 0,
  "I2" + "I3" -> 2*"I3",
  "I3" -> 0
};
rts = {
  la*s,
  mu*s,
  be1*s*i1,
  mu1*i1,
  be2*i1*i2,
  mu2*i2,
  be3*i2*i3,
  mu3*i3
};
\end{verbatim}
Each tier $\mathcal M_j=\{I_j\}$ is exclusive, and the corresponding
restricted submatrix $A_{\mathcal M_j}$ admits a positive flux $v_j$
with $A_{\mathcal M_j}v_j>0$.
The system therefore satisfies both exclusivity and stoichiometric
autocatalysis on each tier, forming a hierarchical SDAS with a tree IGMS.
\end{example}

\begin{example}[Cross-tier catalysis coupling tiers]
\begin{verbatim}
RN = {
  0 -> "S",
  "S" -> 0,
  "S" + "I1" -> 2*"I1",
  "I1" -> 0,
  "S" + "I2" -> 2*"I2",
  "I2" -> 0,
  "I2" + "I1" -> "I2" + "I3",
  "I3" -> 0
};
rts = {
  la*s,
  mu*s,
  be1*s*i1,
  mu1*i1,
  be2*s*i2,
  mu2*i2,
  si*i2*i1,
  mu3*i3
};
\end{verbatim}
The subnetwork is exclusive in
$\mathcal M=\{I_1,I_2,I_3\}$ since all reactions consume at least one
member of~$\mathcal M$.
The cross-tier reaction \(I_2+I_1\to I_2+I_3\) introduces additional
coupling so that $A_{\mathcal M}$ is no longer block-triangular.
Nevertheless, a positive flux $v>0$ with $A_{\mathcal M}v>0$ can still
exist, indicating stoichiometric autocatalysis despite the loss of
tier separation.
\end{example}

\paragraph{Summary}
In these examples, \emph{exclusive autocatalysis} ensures that each
subnetwork ceases to function when its internal species vanish,
whereas \emph{stoichiometric autocatalysis} verifies that a strictly
positive internal flux exists.
The distinction clarifies why tree IGMS can lack SDAS
and why cyclic IGMS can still realize it.

Tree-structured IGMS ensures hierarchical dependency but not autocatalysis.
SDAS requires a strictly positive stoichiometric flux.
Cyclic IGMS may still sustain SDAS through joint autocatalytic cores.

\subsection{The structural point of view on positive ODEs}\lbl{s:DSR}

In this last section, we   make a digression on  the so called ``structural/robust" features of positive ODEs (like the Petri graph), which are fundamental in CRNT, and which we speculate might be useful in ME as well.

For an informal definition, a structural aspect is a feature which is common to all the ODEs in a subclass, and  in  particular, it may be obtained by asking an AI to produce it. For example, the figure \ref{fig:Petri_triangle} below (unrelated to the rest of the paper) was achieved  by typing ``Claude, draw  the Volpert graph of the pasted RN example, arranging the reactions in a triangle, and the species in two nested triangles, such that the outer one forms a regular hexagon together with the  reactions  triangle":

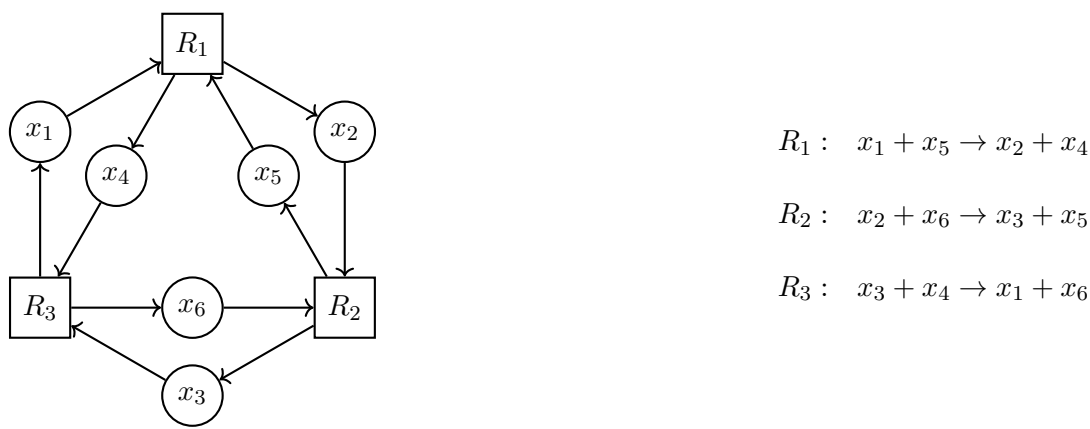
\begin{figure}[h!]
\centering
\begin{minipage}{0.55\textwidth}
\centering
\begin{tikzpicture}[
  scale=0.9,
  place/.style={circle, draw, minimum size=0.8cm, thick},
  transition/.style={rectangle, draw, minimum size=0.8cm, thick},
  arrow/.style={->, thick}
]

\node[transition] (R1) at (0, 2.6) {$R_1$};
\node[place] (x1) at (-2.25, 1.3) {$x_1$};
\node[place] (x2) at (2.25, 1.3) {$x_2$};
\node[transition] (R2) at (2.25, -1.3) {$R_2$};
\node[place] (x3) at (0, -2.6) {$x_3$};
\node[transition] (R3) at (-2.25, -1.3) {$R_3$};

\node[place] (x5) at (1.125, 0.65) {$x_5$};  
\node[place] (x6) at (0, -1.3) {$x_6$};      
\node[place] (x4) at (-1.125, 0.65) {$x_4$}; 

\draw[arrow] (x1) -- (R1);
\draw[arrow] (x5) -- (R1);
\draw[arrow] (R1) -- (x2);
\draw[arrow] (R1) -- (x4);

\draw[arrow] (x2) -- (R2);
\draw[arrow] (x6) -- (R2);
\draw[arrow] (R2) -- (x3);
\draw[arrow] (R2) -- (x5);

\draw[arrow] (x3) -- (R3);
\draw[arrow] (x4) -- (R3);
\draw[arrow] (R3) -- (x1);
\draw[arrow] (R3) -- (x6);

\end{tikzpicture}
\end{minipage}%
\hfill
\begin{minipage}{0.4\textwidth}
\centering
\begin{align*}
R_1: &\quad x_1 + x_5 \to x_2 + x_4 \\[1em]
R_2: &\quad x_2 + x_6 \to x_3 + x_5 \\[1em]
R_3: &\quad x_3 + x_4 \to x_1 + x_6
\end{align*}
\end{minipage}
\caption{Petri net of a CRN. The three reactions/transitions $R_1,R_2,R_3$ appear as squares arranged in a triangle, and the six species/places $x_1,\dots,x_6$ appear as circles.}
\label{fig:Petri_triangle}
\end{figure}

 Now, the ME literature has rarely paid attention to structural features, with the notable exception of the (F,V,K) NGMs, which we consider a structural feature of all ME models.

Historically, ME ODE research is divided in symbolic analyses of small examples for all possible \para s (of which we also bring a sample in this paper), and in numerical analyses of large examples, using ``estimated realistic values" of the constants. CRNT brings two more points of view.  One is identifying subclasses of ODEs for which stability/multistationarity/global attractivity results hold for any values  of  the positive \para s of the model. The second, more recent, is identifying results which hold for any ``admissible rates" (differentiable, monotone, etc) -- see \cite{al2014robust,AAunc,AAnew,AASrob,Ali23gr}. Neither of these particular two directions has yet proven useful  in ME, but in this paper, we offer  a first structural contribution of CRNT to ME, the construction of the DFE as union of minimal siphons, answering a question we asked in \cite{AAlaws}.  As an aside,  note that once the DFE variables are provided, most AI's will have no difficulty
producing the (F,V,K) matrices, for any ME model, specified either in ODE or in RN form. By defining the DFE, and using the algorithmic choice of F in \cite{AABJ},  we have made this computation automatic, as a function  of the ODE.  We  speculate that other CRNT structures, like the Petri graph, will also become useful in future research.


\sec{Conclusion}

In this work, we aimed to  address the fragmented landscape of three traditionally separate disciplines: mathematical epidemiology, chemical reaction network theory, and ecological modeling by proposing a unifying theoretical framework.

This work demonstrates that epidemic strains are fundamentally analogous to minimal siphons, that govern boundary behavior in positive dynamical systems. This correspondence, encoded in the stoichiometric matrices that define how species interact,  influences the block structure of next-generation matrices.

What emerges from this synthesis is a rigorous pathway for analyzing multi-strain epidemic systems through the lens of reproduction functions. These functions, evaluated at various boundary equilibria, capture the invasion potential of competing strains and naturally extend the classical threshold theorem of \cite{Van} to settings where multiple pathogens coexist or compete.
The integration of reaction network formalism into epidemiology also brings algorithmic advantages, illustrated in  the accompanying software package via symbolic and numerical tools.

In conclusion, this work suggests that the tools of CRNT are not incidental conveniences, but fundamental structures governing multi-strain epidemic behavior. By making these tools accessible to epidemiologists and ecologists, and by demonstrating their equivalence to reproduction number based reasoning, the framework presented here contributes to a unified mathematical theory of biological interaction networks.

Looking forward, several directions promise to extend this framework's scope.
 Extending beyond mass-action kinetics to frequency-dependent transmission, saturating incidence, or resource-mediated interactions will broaden applicability to contexts where density-dependent assumptions fail.

\textbf{acknowledgments: } We thank Murad Banaji,  Jackie Cresson, Chuanhou Gao, Vincent Du, Praful Gagrani, Andrei Halanay, Hyukpyo Hong, Matthew Johnston, Sei Zhen Kong, M. Ali Al-Radhawi, Eduardo Sontag,  Nicola Vassena, Polly Yu and Xiaoyu Zhang for useful interactions.

\bibliographystyle{apalike}
\bibliography{ref}
\end{document}